\documentclass[DIV=14,letterpaper]{scrartcl}

\date{}

\usepackage{tikz}
\usetikzlibrary{calc}
\usetikzlibrary{matrix}
\usepackage{amsmath,amssymb,amsfonts,amsthm,subfig}
\usepackage[pdftex]{hyperref} 

\usepackage{graphicx}

\newcommand{\V}[1]{\mbox{\boldmath $ #1 $}}

\newcommand{\bey}{\begin{eqnarray}}
\newcommand{\eey}{\end{eqnarray}}
\newcommand{\nn}{\nonumber}
\newcommand{\beq}{\begin{equation}}
\newcommand{\eeq}{\end{equation}}
\theoremstyle{plain}
\newtheorem{thm}{\hspace{6mm}Theorem}[section]

\newtheorem{lem}{\hspace{6mm}Lemma}[section]

\theoremstyle{definition}
\newtheorem{de}{\hspace{6mm}Definition}[section]
\theoremstyle{remark}
\newtheorem{exam}{\hspace{6mm}Example}[section]
\newtheorem{rem}{\hspace{6mm}Remark}[section]
\newcommand{\proofend}{\mbox{ }\hfill \raisebox{.4ex}{\framebox[1ex]{}}}

\title{Sign-preserving of principal eigenfunctions in P1 finite element approximation
of eigenvalue problems of second-order elliptic operators}

\author{Weizhang~Huang
\thanks{Department of Mathematics, the University of Kansas, Lawrence, KS 66045, U.S.A. 
({\em huang@math.ku.edu}).}
}

\begin{document}
\vskip 1cm
\maketitle

\begin{abstract}
This paper is concerned with the P1 finite element approximation of the eigenvalue problem of
second-order elliptic operators subject to the Dirichlet boundary condition. The focus is on
the preservation of basic properties of the principal eigenvalue and eigenfunctions of
continuous problems. It is shown that when the stiffness matrix is an irreducible $M$-matrix,
the algebraic eigenvalue problem maintains those properties such as the smallest eigenvalue
being real and simple and the corresponding eigenfunctions being either positive or negative
inside the physical domain. Mesh conditions leading to such a stiffness matrix are
also studied. A sufficient condition is that the mesh is simplicial, acute when measured in the metric specified
by the inverse of the diffusion matrix, and interiorly connected. The acute requirement can be
replaced by the Delaunay condition in two dimensions. Numerical results are presented to verify the theoretical findings.
\end{abstract}



\noindent{\bf Abbreviated title.}
Sign-preserving of principal eigenfunctions in FEM

\section{Introduction}

We are concerned with the P1 finite element approximation of the eigenvalue problem of a general second-order
elliptic operator
\beq
\begin{cases}
 L u \equiv -\nabla \cdot (\mathbb{D} \, \nabla u) + \V { b } \cdot \nabla u + c \, u = \lambda u,&\quad \mbox{ in } \Omega \\
 u = 0,&\quad \mbox{ on } \partial\Omega
 \end{cases}
 \label{eigen-1}
\eeq
where $\Omega \subset \mathbb{R}^d$ ($d \ge 1$) is a polyhedron and $\mathbb{D} = \mathbb{D}( \V{x} ): \Omega
\to \mathbb{R}^{d \times d}$, $\V{b} = \V{b}(\V{x}): \Omega \to \mathbb{R}^d$, and $c=c(\V{x}): \Omega \to \mathbb{R}$
are given, sufficiently smooth functions.
We assume that $\mathbb{D}$ is symmetric and strictly positive definite on $\Omega$
and the functions $\V{b}$ and $c$ satisfy
\beq
c(\V{x}) - \frac{1}{2} \nabla \cdot \V{b}(\V{x}) \geq 0 ,\qquad \forall \V{x}\in \Omega .
\label{eigen-2}
\eeq
Note that the condition (\ref{eigen-2}) is not essential. We can always make them satisfied
by adding a large positive number to the function $c(\V{x})$. The original and shifted problems will have
the same eigenfunctions and the eigenvalues of the former can be obtained by shifting the eigenvalues of the latter.

The eigenvalue problem (\ref{eigen-1}) is not self-adjoint in general. Nevertheless, it is known
(e.g., see Lemma~\ref{lem:pde-eigen} below or
Evans \cite[Theorem 2 on Page 336 and Theorem 3 on page 340]{Eva98})
that the principal eigenvalue (that is, the smallest
eigenvalue in modulus) is real and simple and the principal eigenfunctions
(that is, the eigenfunctions corresponding to the principal
eigenvalue) are either positive or negative in $\Omega$.
Since the principal eigenvalues typically represent the ground state of a physical system or correspond to
the most unstable mode in stability or sensitivity analysis, it is of practical and theoretical importance to
study when a numerical approximation preserves these properties of the principal eigenvalue and eigenfunctions
and especially the sign of the principal eigenfunctions.

The objective of this paper is to investigate when a P1 finite element approximation of (\ref{eigen-1})
on a simplicial mesh preserves the basic properties of the principal eigenvalue and eigenfunctions of the continuous
problem. We shall show that most of the basic properties of the principal eigenvalue and eigenfunctions 
are preserved in the P1 finite element approximation provided that the resulting stiffness matrix is an irreducible
$M$-matrix (cf. Theorem~\ref{thm:fem-eigen}). Particularly, the principal eigenvalue of the discrete system
is real and algebraically and geometrically simple and the corresponding eigenfunctions are either positive
or negative (sign-preserving) in the physical domain. 
Several sufficient mesh conditions are proposed (Theorem~\ref{thm:irreducibleM}) for the stiffness matrix to be
an irreducible $M$-matrix.

We point out that there is a vast literature on the finite element approximation of
differential eigenvalue problems and most of it is on convergence analysis;
e.g., see Babu{\v{s}}ka and Osborn \cite{BO1991},
Boffi \cite{Bof2010}, and Boffi et al. \cite{Boffi2010b}, and references therein.
Early work includes Birkhoff et al. \cite{BdBSW1966}, Fix \cite{Fix1973}, 
and Babu{\v{s}}ka and Osborn \cite{BO1989}.
We also point out some interesting recent work
\cite{DXZ2008,DZ2008,HuJun2012,LLX2012,MM2011,NZ2012,XZ2001,YB2011,YZL2010}.

An outline of the paper is as follows. The P1 finite element approximation of (\ref{eigen-1}) is presented in \S\ref{SEC:p1fem}
and the preservation of the basic properties of the principal eigenvalue and eigenfunctions in the P1 finite element
approximation is studied in \S\ref{SEC:preservation}. \S\ref{SEC:meshcond} is devoted to the study of 
mesh conditions under which the stiffness matrix is ensured to be an irreducible $M$-matrix, followed by
numerical examples in \S\ref{SEC:numerics}. The conclusions are drawn in \S\ref{SEC:conclusions}.

\section{P1 finite element formulation}
\label{SEC:p1fem}

The weak formulation of the eigenvalue problem (\ref{eigen-1}) is to find $\lambda \in \mathbb{C}$ and nonzero
(and possibly complex) function $u\in H_0^1(\Omega)$ such that
\beq
(\mathbb{D} \nabla u, \nabla v) + (\V { b } \cdot \nabla u, v) + (c\, u, v) = \lambda (u, v),\quad \forall v \in H_0^1(\Omega)
\label{eigen-3}
\eeq
where $(\cdot, \cdot)$ denotes the $L^2$ inner product.
For the P1 finite element approximation, we assume that
an affine family of simplicial mesh $\{ \mathcal{T}_h \}$ is given for $\Omega$.
Denote by $V^h\subset H_0^1(\Omega)$ the standard P1 finite element space associated
with a mesh $\mathcal{T}_h$. A P1 finite element approximation to the eigenvalue problem (\ref{eigen-3}) is to
find $\lambda^h \in \mathbb{C}$ and nonzero (and possibly complex) function $u^h \in V^h$ such that
\beq
(\mathbb{D} \nabla u^h, \nabla v^h) + (\V { b } \cdot \nabla u^h, v^h) + (c\, u^h, v^h) = \lambda^h (u^h, v^h),
\quad \forall \; v^h \in V^h .
\label{fem-1}
\eeq

Scheme (\ref{fem-1}) can be expressed in a matrix form.
Denote the numbers of the elements and the interior vertices of $\mathcal{T}_h$
by $N$ and $N_{v}$, respectively. Assume that the vertices are ordered in such a way that
the first $N_{v}$ vertices are the interior ones. Then, $V^h$ and $u^h$ can be expressed as
\[
V^h = \text{span} \{ \phi_1, \cdots, \phi_{N_{v}} \} ,\quad u^h = \sum_{k=1}^{N_{v}} u_k \phi_k  , 
\]
where $\phi_k$ denotes the P1 basis function associated with the $k^{\text{th}}$ vertex.
Substituting the above expression into (\ref{fem-1}) and taking $v^h = \phi_j$ ($j=1, ..., N_{v}$),
we obtain the algebraic eigenvalue problem
\beq
A \, \V{u} = \lambda^h B \V{u},
\label{fem-2}
\eeq
where $\V{u} = (u_1,..., u_{N_{v}})^T$, the stiffness matrix $A$ and the mass matrix $B$ are given by
\begin{align}
& A_{j k}  = \sum_{K \in \mathcal{T}_h} \int_K d \V{x} \left ( \frac{}{} (\nabla \phi_j)^{T} \, \mathbb{D}_K \,
\nabla \phi_k +  \phi_j \; (\V{b} \cdot \nabla \phi_k) + c \; \phi_j \, \phi_k \right ) ,
\quad j, k =1, ..., N_{v}
\label{A-1}
\\
& B_{j k} = \sum_{K \in \mathcal{T}_h} \int_{K} \phi_j \phi_k d \V{x},\quad j, k =1, ..., N_{v} 
\label{B-1}
\end{align}
and $\mathbb{D}_K$ is the average of $\mathbb{D}$ over $K$, i.e.,
\beq
\mathbb{D}_K = \frac{1}{|K|} \int_K \mathbb{D}(\V{x})\, d \V{x} .
\label{D-1}
\eeq

The convergence of finite element approximation of (\ref{eigen-3}) has been extensively studied
(e.g., see \cite{Bof2010}). We can expect that the principal eigenvalue of (\ref{fem-1})
converges to that of the continuous problem (\ref{eigen-3}) at $\mathcal{O}(N^{-\frac{2}{d}})$ (second order)
as $N \to \infty$. On the other hand, the preservation of the basic properties
of the principal eigenvalue and eigenfunctions by the discrete system has not been studied so far (to our best knowledge). Our goal is to establish conditions (on the stiffness matrix and the mesh) under which
the discrete eigenvalue problem (\ref{fem-2})
preserves those properties.

\section{Preservation of basic properties of the principal eigenvalue and eigenfunctions}
\label{SEC:preservation}

In this section we describe basic properties of the principal eigenvalue and eigenfunctions of
the continuous problem (\ref{eigen-3}) (Lemma~\ref{lem:pde-eigen}) and show (Theorem~\ref{thm:fem-eigen})
that those properties are preserved by the discrete eigenvalue problem (\ref{fem-1})
provided that the stiffness matrix $A$ is an irreducible $M$-matrix.
The mesh conditions to ensure an irreducible $M$-matrix stiffness matrix for the P1 finite element approximation
will be studied in the next section.

\begin{lem}
\label{lem:pde-eigen}
The principal eigenvalue $\lambda_1 $ and the corresponding eigenfunctions of
the eigenvalue problem (\ref{eigen-3}) have the following properties.
\begin{itemize}
\item[(a)] $\lambda_1$ is real;
\item[(b)] There is an eigenfunction $u_1 \in H_0^1(\Omega)$ associated with $\lambda_1$
	with $u_1(\V{x}) > 0$ for all $\V{x} \in \Omega$;
\item[(c)] $\lambda_1$ is simple, that is, if $u$ is an eigenfunction associated with $\lambda_1$,
	then $u$ is a multiple of $u_1$;
\item[(d)] $\lambda_1 = F (u_1)  >  0$, where $F(\cdot)$ is defined as
\beq
F(v) = \frac{(\mathbb{D} \nabla v, \nabla v) + ((c-\frac{1}{2} \nabla \cdot \V{b}) v, v)}{(v,v)}  ;
\label{v-3}
\eeq
\item[(e)] $ \text{Re}(\lambda)\ge \lambda_1$ for every eigenvalue $\lambda$;
\item[(f)] For the symmetric situation (with $\V{b} = 0$), there holds the variational principle
	\beq
	\lambda_1 = \min\limits_{\substack{v \in H_0^1(\Omega)\\ v \not\equiv 0,\; \text{real} }} F (v) .
	\label{lem:pde-eigen-1}
	\eeq
\end{itemize}
\end{lem}

{\bf Proof.} (a), (b), (c), (e), and (f) are the standard results for second-order elliptic operators; e.g, see
\cite[Theorem 2 on Page 336 and Theorem 3 on page 340]{Eva98}. (d) follows from equation (\ref{eigen-3})
(with $u = v = u_1$), integration by parts,
the assumption (\ref{eigen-2}), Poincar\'{e}'s inequality, and the fact that $\lambda_1$
and $u_1$ are real.
\proofend

\vspace{10pt}

\begin{thm}
\label{thm:fem-eigen}
For the finite element eigenvalue problem (\ref{fem-1}), if the stiffness matrix $A$ is an irreducible $M$-matrix, then
the principal eigenvalue $\lambda_1^h$ and the corresponding eigenfunctions have the following properties.
\begin{itemize}
\item[(a)] $\lambda_1^h$ is real;
\item[(b)] There is an eigenfunction $u_1^h \in V^h$ associated with $\lambda_1^h$
	with $u_1^h(\V{x}) > 0$ for all $\V{x} \in \Omega$;
\item[(c)] $\lambda_1^h$ is algebraically (and geometrically) simple,
	that is, if $u^h$ is an eigenfunction associated with $\lambda_1^h$,
	then $u^h$ is a multiple of $u_1^h$;
\item[(d)] $\lambda_1^h = F (u_1^h)  >  0$, where $F(\cdot)$ is defined in (\ref{v-3});
\item[(e)] $\text{Re}(\lambda^h) > 0$ and $ |\lambda^h|\ge \lambda_1$ for every eigenvalue $\lambda^h$;
\item[(f)] For the symmetric situation (with $\V{b} = 0$), there holds the variational principle
	\beq
	\lambda_1^h = \min\limits_{\substack{v^h \in V^h\\ v^h \not\equiv 0,\; \text{real} }} F (v^h) .
	\label{thm:fem-eigen-1}
	\eeq
\end{itemize}
\end{thm}

{\bf Proof.}
The finite element eigenvalue problem (\ref{fem-1}) or (\ref{fem-2}) is mathematically equivalent to
\beq
A^{-1} B \V{u} = \frac{1}{\lambda^h} \V{u} .
\label{thm:fem-eigen-2}
\eeq
Since $A$ is an irreducible $M$-matrix, $A^{-1}$ is positive, i.e., $A^{-1} > 0$ (in the elementwise sense).
From (\ref{B-1}), it is obvious that each column of $B$ has at least one non-zero entry.
Thus, we have $A^{-1} B > 0$. We also notice that $u^h(\V{x}) = \sum_{j=1}^{N_{v}} u_j \phi_j(\V{x}) > 0$
for all $\V{x} \in \Omega$ if and only if $\V{u} = (u_1, ..., u_{N_v})^T > 0$.
Then, (a), (b), and (c) follow from Perron's Theorem (for positive matrices; e.g., see \cite[8.2.11]{HJ1985}).
(d) follows from the equation (\ref{fem-1}) (with $u^h = v^h = u_1^h$), integration by parts,
the assumption (\ref{eigen-2}), Poincar\'{e}'s inequality, and the fact
that $\lambda_1^h$ and $u_1^h$ are real.
For (e), the property $\text{Re}(\lambda^h) > 0$ is a consequence of the fact that $A$ is an $M$-matrix and
$B$ is a symmetric and positive definite matrix. The other property follows from Perron's Theorem.

Next, we show that (f) holds. For this case, $A$ is symmetric. From Perron's Theorem, the eigenvalues of
(\ref{fem-2}) can be ordered as
\[
0 < \lambda_1^h < \lambda_2^h \le \cdots \le \lambda_{N_v}^h .
\]
Denote the corresponding normalized eigenvectors by $u_{j}^h$ (or $\V{u}_j$ in vector form) ($j = 1, ..., N_v$).
Notice that they satisfy $\V{u}_j^T B \V{u}_k = \delta_{jk}$. Then, any function $v^h \in V^h$ can be expressed into
\[
v^h = \sum_{j=1}^{N_v} d_j u_j^h \quad \text{ or }\quad 
\V{v} = \sum_{j=1}^{N_v} d_j \V{u}_j .
\]
From the orthogonality of the eigenfunctions, we have
\[
F(v^h) = \frac{\V{v}^T A \V{v}}{\V{v}^T B \V{v}} = \frac{\sum_j d_j^2 \lambda_j^h}{\sum_j d_j^2}
\ge \lambda_1^h .
\]
Combining this with (d) gives (\ref{thm:fem-eigen-1}).
\proofend

\vspace{10pt}

\begin{rem}
\label{rem3.1}
We note that the properties in Theorem~\ref{thm:fem-eigen}(e) are different from the stronger property
$\text{Re}(\lambda^h) \ge \lambda_1^h$ (cf. Lemma~\ref{lem:pde-eigen}(e)).
There are two cases we can show that the discrete system has the latter property when $A$ is
an irreducible $M$-matrix.
The first case is the symmetric case. In this case, $\text{Re}(\lambda^h) = \lambda^h$, and
$\text{Re}(\lambda^h) > 0$ and $ |\lambda^h|\ge \lambda_1$
imply $\text{Re}(\lambda^h) \ge \lambda_1^h$. The other case is to use the lumped mass matrix
(denoted by $\tilde{B}$) instead of the full mass matrix $B$. Since $\tilde{B}$ is diagonal,
$\tilde{B}^{-1} A$ is also an irreducible $M$-matrix, which implies $\text{Re}(\lambda^h) \ge \lambda_1^h$
(e.g., see Elhashash and Szyld \cite[Theorem 3.1]{ES2008}).

For the general nonsymmetric situation, we are unable to show that the discrete system has the property
$\text{Re}(\lambda^h) \ge \lambda_1^h$ for every eigenvalue $\lambda^h$
although our limited numerical experiment shows that the system does satisfy the property (cf. Fig.~\ref{ex5.2-f2}).
\proofend
\end{rem}

\vspace{10pt}

Theorem~\ref{thm:fem-eigen} states that if $A$ is an irreducible $M$-matrix, then the P1 finite element approximation
(\ref{fem-1}) essentially retains most of the properties listed in Lemma~\ref{lem:pde-eigen} for the principal eigenvalue
and eigenfunctions. In the next section we study the mesh conditions
to ensure that the P1 finite element stiffness matrix be an irreducible $M$-matrix. 

\section{Mesh conditions for irreducible $M$-matrix stiffness matrix}
\label{SEC:meshcond}

We first study mesh conditions to ensure $A$ to be an $M$-matrix. This issue is closely related to
the preservation of the maximum principle for boundary value problems. The latter has been studied extensively
in the past; for example, see \cite{BKK08,BE04,Cia70,CR73,KK09,KKK07,KL95,Let92,Sto86,SF73,WaZh11,XZ99}
for isotropic diffusion problems ($\mathbb{D} = \alpha(\V{x}) I$ with $\alpha(\V{x})$ being a scalar function) and
\cite{DDS04,GL09,GYK05,Hua10,KSS09,LePot09,LH10,LSS07,LSSV07,LS08,LHQ2012,MD06,SH07,
ShYu2011,YuSh2008,ZZS2013} for anisotropic diffusion problems.

In the following we quote a result from Lu et al. \cite{LHQ2012}.  We first introduce some notation.
For any simplicial element $K$, we denote the inner normal to face $S_j^K$
(the face not containing the $j^{\text{th}}$ vertex of $K$) by $\V{q}_j^K$. The dihedral angle
in the metric $\mathbb{D}^{-1}$ between faces $S_j^K$ and $S_k^K$ ($j \neq k$) can be computed as
\beq
\alpha_{j k, \mathbb{D}^{-1}}^K = - \frac{(\V{q}_j^K)^T \mathbb{D}_K \V{q}_k^K}
{\sqrt{(\V{q}_j^K)^T \mathbb{D}_K \V{q}_j^K \cdot (\V{q}_k^K)^T \mathbb{D}_K \V{q}_k^K} } .
\label{alpha-1}
\eeq
The maximum dihedral angle in the metric $\mathbb{D}^{-1}$ for $K$ is defined as
\beq
\alpha_{\max, \mathbb{D}^{-1}}^K = \max\limits_{j,k=1, ..., d+1, j \neq k} \alpha_{j k, \mathbb{D}^{-1}}^K .
\label{alpha-2}
\eeq
The diameter (i.e., the largest edge length in the Euclidean metric) of $K$ is denoted by $h_K$.

\begin{lem}
\label{lem:LHQ12}
If the mesh satisfies
\beq
0 < \alpha_{\max, \mathbb{D}^{-1}}^K \le \arccos\left (
 \frac{h_K}{\lambda_{min}(\mathbb{D}_K)} \cdot \frac{\|\V{b}\|_{L^{\infty}(K)} }{(d+1)}
+ \frac{h_K^2}{\lambda_{min}(\mathbb{D}_K)}\cdot \frac{\|c\|_{L^{\infty}(K)} }{(d+1)(d+2)}
\right ), \quad \forall K \in \mathcal{T}_h
\label{lem:LHQ12-1}
\eeq
then, the stiffness matrix $A$ is an $M$-matrix.

In 2D, the above condition can be replaced by a Delaunay-type condition
\bey
&& 0 < \frac{1}{2} \left [ \frac{}{} \alpha_{jk, \mathbb{D}^{-1}}^K + \alpha_{jk, \mathbb{D}^{-1}}^{K'}
\right . \nn \\
&& \quad
+ \text{ arccot}\left (\sqrt{\frac{\mbox{ det}(\mathbb{D}_{K'})}{\mbox{ det}(\mathbb{D}_{K})}}
\cot (\alpha_{jk, \mathbb{D}^{-1}}^{K'} ) - \frac{2\; \Theta(K,K') }{\sqrt{\mbox{ det}(\mathbb{D}_{K})}} \right )
\nn \\
&& \quad \left . + \text{ arccot}\left (\sqrt{\frac{\mbox{ det}(\mathbb{D}_{K})}{\mbox{ det}(\mathbb{D}_{K'})}}
\cot (\alpha_{jk, \mathbb{D}^{-1}}^K ) - \frac{2\; \Theta(K,K') }{\sqrt{\mbox{ det}(\mathbb{D}_{K'})}} \right )
\frac{}{} \right ] \le \pi 
\label{lem:LHQ12-2}
\eey
for every internal edge $e_{jk}$ connecting the $j^{\text{th}}$ and $k^{\text{th}}$ vertices. Here, $K$
and $K'$ are the elements sharing the common edge $e_{jk}$, $\alpha_{jk, \mathbb{D}^{-1}}^K$
and $\alpha_{jk, \mathbb{D}^{-1}}^{K'}$ are the angles in $K$ and $K'$ that face the edge, and
\beq
\Theta(K,K') = \frac{h_K \; \|\V{b}\|_{L^{\infty}(K)}}{ (d+1)} + \frac{h_K^2 \; \|c\|_{L^{\infty}(K)}}{(d+1)(d+2)}
+ \frac{h_{K'} \; \|\V{b}\|_{L^{\infty}(K)}}{ (d+1)} + \; \frac{h_{K'}^2 \; \|c\|_{L^{\infty}(K)}}{(d+1)(d+2)} .
\label{lem:LHQ12-3}
\eeq
\end{lem}

{\bf Proof.} This result was proven in Lu et al. \cite[Theorems 1 and 2]{LHQ2012}. For completeness,
we give a proof here. The proof is also useful in the study of irreducibility of the stiffness matrix,
see Theorem~\ref{thm:irreducibleM}.

We first show that $A$ is a Z-matrix; i.e.,
\bey
a_{j k} &\leq& 0, \;\; \forall\; j \neq k, \; j, k = 1, ..., N_{v}
\nn
\\
a_{jj} &\geq& 0, \;\; j = 1,..., N_{v} .
\nn
\eey
Recall from Ciarlet \cite[Page 201]{Cia78} that
\beq
\int\limits_{K \in \omega_j} \phi_j d \V{x} = \frac{|K|}{d+1},\quad
\int\limits_{K \in \omega_j \cap \omega_k} \phi_j \phi_k d \V{x} = \frac{|K|}{(d+1)(d+2)} ,\quad j \ne k
\label{ciarlet78-1}
\eeq
where $\omega_j$ and $\omega_k$ are the element patches associated with the $j^{\text{th}}$ and $k^{\text{th}}$ vertices,
respectively. For $j\ne k$, from (\ref{A-1}) we have
\[
a_{jk} = \sum_{K \in \omega_j \cap \omega_k} \left (
|K| \;(\nabla \phi_j)^T \; \mathbb{D}_K \; \nabla \phi_k
+ \int_K \; \phi_j \; (\V{b} \cdot \nabla \phi_k ) d \V{x}
+ \; \int_K
\; c \; \phi_j \; \phi_k d \V{x} \right ) .
\]
From \cite[Lemmas 1 and 3]{LHQ2012},
\[
\left. \nabla \phi_j \right |_K = - \frac{1}{h_j^K} \frac{\V{q}_j^K}{\sqrt{(\V{q}_j^K)^T \V{q}_j^K }},\quad 
\left. (\nabla \phi_j)^T \; \mathbb{D}_K \; \nabla \phi_k\right |_K =
- \frac{\cos(\alpha_{jk,\mathbb{D}^{-1}}^K)}{h_{j, \mathbb{D}^{-1}}^K h_{k, \mathbb{D}^{-1}}^K}, 
\]
where $h_j^K$ and $h_{j, \mathbb{D}^{-1}}^K$ are the $j^{\text{th}}$ altitude of $K$ in the Euclidean metric and
the metric specified by $\mathbb{D}^{-1}$, respectively. They are related by
\[
\frac{h_j^K}{\sqrt{\lambda_{\max} (\mathbb{D}_K)}}  \le 
h_{j, \mathbb{D}^{-1}}^K \le \frac{h_j^K}{\sqrt{\lambda_{\min} (\mathbb{D}_K)}} .
\]
Combining the above results, we have
\begin{align}
a_{jk} &\leq \sum_{K \in \omega_j \cap \omega_k} \left ( - \frac{|K|}{h_{j, \mathbb{D}^{-1}}^K h_{k, \mathbb{D}^{-1}}^K}
\cos(\alpha_{jk,\mathbb{D}^{-1}}^K)
+  \frac{\|\V{b}\|_{L^{\infty}(K)}}{h_k^K} \int_K \phi_j   d \V{x}
+ \|c\|_{L^{\infty}(K)} \int_K \; \phi_j \; \phi_k \; d \V{x} \right )
\notag \\
&= \sum_{K \in \omega_j \cap \omega_k} \left ( - \frac{|K|}{h_{j, \mathbb{D}^{-1}}^K h_{k, \mathbb{D}^{-1}}^K}
\cos(\alpha_{jk,\mathbb{D}^{-1}}^K)
+  \frac{|K|\; \|\V{b}\|_{L^{\infty}(K)}}{(d+1) h_k^K}
+ \frac{|K| \; \|c\|_{L^{\infty}(K)}}{(d+1)(d+2)}  \right )
\notag \\
&= \sum_{K \in \omega_j \cap \omega_k} \frac{|K|}{h_{j, \mathbb{D}^{-1}}^K h_{k, \mathbb{D}^{-1}}^K} 
\left ( -  \cos(\alpha_{jk,\mathbb{D}^{-1}}^K)
+  \frac{h_{j, \mathbb{D}^{-1}}^K h_{k, \mathbb{D}^{-1}}^K \|\V{b}\|_{L^{\infty}(K)}}{(d+1) h_k^K}
+ \frac{h_{j, \mathbb{D}^{-1}}^K h_{k, \mathbb{D}^{-1}}^K \|c\|_{L^{\infty}(K)}}{(d+1)(d+2)}  \right )
\notag \\
&\leq \sum_{K \in \omega_j \cap \omega_k} \frac{|K|}{h_{j, \mathbb{D}^{-1}}^K h_{k, \mathbb{D}^{-1}}^K} 
\left ( -  \cos(\alpha_{\max,\mathbb{D}^{-1}}^K)
+  \frac{h_K \|\V{b}\|_{L^{\infty}(K)}}{(d+1) \lambda_{\min} (\mathbb{D}_K) }
+ \frac{h_K^2 \|c\|_{L^{\infty}(K)}}{(d+1)(d+2) \lambda_{\min} (\mathbb{D}_K)}  \right ) .
\label{ajk-1}
\end{align}
Thus, $a_{j,k} \le 0$ when (\ref{lem:LHQ12-1}) is satisfied.

In two dimensions, notice that there are only two elements in $\omega_j \cap \omega_k$ which share
the common edge $e_{jk}$. Denote these elements by $K$ and $K'$. Similarly, we can get
\begin{align}
a_{jk} & \le - \frac{\mbox{det}(\mathbb{D}_K)^{\frac 1 2}}{2} \cot(\alpha_{jk, \mathbb{D}^{-1}}^K)
- \frac{\mbox{det}(\mathbb{D}_{K'})^{\frac 1 2}}{2} \cot(\alpha_{jk, \mathbb{D}^{-1}}^{K'})
\nn \\
& \quad + \frac{h_K \; \|\V{b}\|_{L^{\infty}(K)}}{ (d+1)} + \frac{h_K^2 \; \|c\|_{L^{\infty}(K)}}{(d+1)(d+2)}
+ \frac{h_{K'} \; \|\V{b}\|_{L^{\infty}(K')}}{ (d+1)} + \frac{h_{K'}^2 \; \|c\|_{L^{\infty}(K')}}{(d+1)(d+2)} .
\label{ajk-2}
\end{align}
It can be shown (e.g., see \cite{Hua10}) that $a_{jk} \le 0$ when (\ref{lem:LHQ12-2}) is satisfied.

For the diagonal entries, we have
\begin{align*}
a_{jj} &= \sum_{K \in \mathcal{T}_h} |K| \;(\nabla \phi_i)^T \; \mathbb{D}_K \; \nabla \phi_i
+ \int_\Omega \; \phi_i \; (\V{b} \cdot \nabla \phi_i ) d \V{x}
+ \int_\Omega\; c \; \phi_i^2 d \V{x}
\\
& \ge \int_{\Omega} \phi_i ( \V{b} \cdot \nabla \phi_i ) d \V{x} + \int_{\Omega} c \; \phi_i^2 d \V{x}
=  \int_{\Omega} ( c - \frac{1}{2} \nabla \cdot \V{b} ) \phi_i^2 d \V{x} \ge 0 .
\end{align*}
Thus, $A$ is a Z-matrix.

We now show that $A$ is an M-matrix by showing that $A$ is positive definite.
For any vector $\V{v}=(v_1, v_2, ..., v_{N_{v}})^T$,
we define $v^h = \sum _{i=1}^{N_{v}}v_i \phi_i \in V^h$. Notice that $\nabla v^h$ is constant on $K$.
As in the proof for $a_{jj}\ge 0$, from (\ref{A-1}) we have
\bey
\V{v}^T A \V{v} &=& \sum_{K \in \mathcal{T}_h} |K| \;( \nabla v^h )^T \; \mathbb{D}_K \; \nabla v^h
+ \int_\Omega v^h \; (\V{b} \cdot \nabla v^h )d \V{x}
+ \int_\Omega c \; (v^h)^2 d \V{x}
\nn \\
&=& \int_\Omega \;( \nabla v^h )^T \; \mathbb{D} \; \nabla v^h d \V{x} +
\int _{\Omega} ( c - \frac{1}{2} \nabla \cdot \V{b} ) ( v^h )^2 d \V{x} \ge 0.
\nn
\eey
Moreover, from the above inequality it is easy to see that $\V{v}^T A \V{v} = 0$ implies 
\[
\int_\Omega \;( \nabla v^h )^T \; \mathbb{D} \; \nabla v^h d \V{x} = 0 .
\]
We thus have $\nabla v^h = 0$ or $v^h = $ constant, which in turn
implies $v^h = 0$ due to the fact that $v^h$ vanishes on $\partial \Omega$.
Hence, $A$ is positive definite.
\proofend

\vspace{10pt}

\begin{rem}
\label{rem4.1}
Loosely speaking, the mesh conditions (\ref{lem:LHQ12-1}) and (\ref{lem:LHQ12-2})
can be written as
\begin{align}
 0 <  \, & \alpha_{\max, \mathbb{D}^{-1}}^K \le \frac{\pi}{2} - C_1 \|\V{b}\|_{L^{\infty}(\Omega)} h
- C_2 \|c\|_{L^{\infty}(\Omega)} h^2,\quad \forall K \in \mathcal{T}_h 
\label{lem:LHQ12-4}
\\
0 < \, & \frac{1}{2} \left [ \frac{}{} \alpha_{jk, \mathbb{D}^{-1}}^K + \alpha_{jk, \mathbb{D}^{-1}}^{K'}
+ \text{ arccot}\left (\sqrt{\frac{\mbox{ det}(\mathbb{D}_{K'})}{\mbox{ det}(\mathbb{D}_{K})}}
\cot (\alpha_{jk, \mathbb{D}^{-1}}^{K'} ) \right ) \right.
\nn \\
& \left .
 + \text{ arccot}\left (\sqrt{\frac{\mbox{ det}(\mathbb{D}_{K})}{\mbox{ det}(\mathbb{D}_{K'})}}
\cot (\alpha_{jk, \mathbb{D}^{-1}}^K ) \right ) \frac{}{} \right ] 
\nn \\
& \qquad \qquad 
\le \pi - C_3 \|\V{b}\|_{L^{\infty}(\Omega)} h
- C_4 \|c\|_{L^{\infty}(\Omega)} h^2, \quad \forall \text{ interior edge $e_{jk}$}
\label{lem:LHQ12-5}
\end{align}
for some positive constants $C_1$, $C_2$, $C_3$, and $C_4$. When
$\mathbb{D} = I$ and $\V{b} \equiv 0$ and $c \equiv 0$, (\ref{lem:LHQ12-5}) becomes
the Delaunay condition, i.e., $0 < \alpha_{jk}^K + \alpha_{jk}^{K'} \le \pi$.
\proofend
\end{rem}

\vspace{10pt}

\begin{rem}
\label{rem4.2}
The conditions (\ref{lem:LHQ12-1}) and (\ref{lem:LHQ12-2}) have several existing
mesh conditions as special examples. They reduce to the mesh conditions of 
Ciarlet and Raviart \cite{CR73} (the nonobtuse angle condition) for isotropic diffusion problems,
Strang and Fix \cite{SF73} (the Delaunay condition) for 2D isotropic diffusion problems, 
Wang and Zhang \cite{WaZh11} for  isotropic diffusion problems with convection and reaction terms, 
Li and Huang \cite{LH10} (the anisotropic nonobtuse angle condition) for anisotropic diffusion problems, and
Huang \cite{Hua10} (a Delaunay-type condition) for 2D anisotropic diffusion problems.
\proofend
\end{rem}

\vspace{10pt}

We now study the irreducibility of the stiffness matrix $A$ using the notion of directed graphs
(e.g., see Berman and Plemmons \cite{BP94}).
The directed graph (denoted by $G(A)$)
of $A$ is defined as a graph consisting of $N_v$ vertices
$P_1, ..., P_{N_v}$, where an edge leads from $P_j$ to $P_k$ if and only if $a_{jk} \ne 0$.
$G(A)$ is said to be strongly connected if 
for any ordered pair $(P_j, P_k)$ of vertices of $G(A)$, there is a sequence of edges which leads from
$P_j$ to $P_k$. Note that in the current situation, the vertices in $G(A)$ have a one-to-one correspondence
to the interior vertices of the mesh.

\begin{de}
\label{interiorly-connected} A mesh is called to be {\it interiorly connected} 
if any two interior vertices of the mesh are connected by a sequence of interior edges.
\proofend
\end{de}

\begin{thm}
\label{thm:irreducibleM}
The stiffness matrix for the $P1$ finite element approximation of (\ref{eigen-3}) is an irreducible $M$-matrix
if the mesh is interiorly connected  and satisfies 
\beq
0 < \alpha_{\max, \mathbb{D}^{-1}}^K < \arccos\left (
 \frac{h_K}{\lambda_{min}(\mathbb{D}_K)} \cdot \frac{\|\V{b}\|_{L^{\infty}(K)} }{(d+1)}
+ \frac{h_K^2}{\lambda_{min}(\mathbb{D}_K)}\cdot \frac{\|c\|_{L^{\infty}(K)} }{(d+1)(d+2)}
\right ), \quad \forall K \in \mathcal{T}_h .
\label{thm:irreducibleM-1}
\eeq

In 2D, the condition (\ref{thm:irreducibleM-1}) can be replaced by a Delaunay-type condition
\bey
&& 0 < \frac{1}{2} \left [ \frac{}{} \alpha_{jk, \mathbb{D}^{-1}}^K + \alpha_{jk, \mathbb{D}^{-1}}^{K'}
\right . \nn \\
&& \quad
+ \text{ arccot}\left (\sqrt{\frac{\mbox{ det}(\mathbb{D}_{K'})}{\mbox{ det}(\mathbb{D}_{K})}}
\cot (\alpha_{jk, \mathbb{D}^{-1}}^{K'} ) - \frac{2\; \Theta(K,K') }{\sqrt{\mbox{ det}(\mathbb{D}_{K})}} \right )
\nn \\
&& \quad \left . + \text{ arccot}\left (\sqrt{\frac{\mbox{ det}(\mathbb{D}_{K})}{\mbox{ det}(\mathbb{D}_{K'})}}
\cot (\alpha_{jk, \mathbb{D}^{-1}}^K ) - \frac{2\; \Theta(K,K') }{\sqrt{\mbox{ det}(\mathbb{D}_{K'})}} \right )
\frac{}{} \right ] < \pi
\label{thm:irreducibleM-2}
\eey
for every internal edge $e_{jk}$ connecting the $j^{\text{th}}$ and $k^{\text{th}}$ vertices, where $K$
and $K'$ are the elements sharing the common edge $e_{jk}$, $\alpha_{jk, \mathbb{D}^{-1}}^K$
and $\alpha_{jk, \mathbb{D}^{-1}}^{K'}$ are the angles in $K$ and $K'$ that face the edge, and
$\Theta(K,K')$ is defined in (\ref{lem:LHQ12-3}).
\end{thm}

{\bf Proof.}
For any pair $(j,k)$ of neighboring mesh vertices, $\omega_j \cap \omega_k \neq \emptyset$. 
From (\ref{ajk-1}), we can see that if (\ref{lem:LHQ12-1}) holds strictly (i.e., (\ref{thm:irreducibleM-1}) holds),
then $a_{jk} < 0$ and $a_{kj} < 0$, that is, $P_j$ and $P_k$ are connected in both directions.
Consequently, if any two vertices of the mesh are connected
by a sequence of interior edges, then $G(A)$ is strongly connected, which in turn implies that
$A$ is irreducible (e.g., see Berman and Plemmons \cite[Theorem~(2.7)]{BP94}).
Combining this and Lemma~\ref{lem:LHQ12}, we have proven that $A$ is an irreducible $M$-matrix.
\proofend

\vspace{20pt}

We now comment on how to generate meshes satisfying (\ref{thm:irreducibleM-1}) or (\ref{thm:irreducibleM-2}).
Since meshes satisfying these conditions are $\mathcal{O}( h \| \V{b} \|_{L^\infty(\Omega)} + h^2 \| c \|_{L^\infty(\Omega)})$
perturbations of acute meshes (or Delaunay meshes in 2D) in the metric $\mathbb{D}^{-1}$ (cf. Remark~\ref{rem4.1}),
 we focus our discussion on the generation of the latter.

When $\mathbb{D}^{-1} = I$, acute or Delaunay meshes in the metric $\mathbb{D}^{-1}$ are simply
acute or Delaunay meshes in the Euclidean metric. Delaunay meshes in 2D
can be generated using many algorithms, e.g., see de Berg et al. \cite{BKOS00}. Moreover,
2D polygonal and 3D polyhedral domains can be partitioned into simplices with acute angles;
e.g., see \cite{BCER95,BEG94,BKKS09,ESU04}.

On the other hand, it is theoretically unknown whether or not acute or Delaunay meshes in a given metric
$\mathbb{D}^{-1} \neq I$ can be generated for general polygonal or polyhedral domains.
Nevertheless, their approximations can be obtained in practice using the notion of (simplicial) $M$-uniform meshes
or uniform meshes in the metric tensor specified by a tensor $M =M(\V{x})$. ($M = \mathbb{D}^{-1}$ for the current
situation.) It is known \cite{Hua06,HR11} that an $M$-uniform mesh satisfies the so-called
equidistribution and alignment conditions
\begin{align}
|K| \text{det}(M_K)^{\frac 1 2} = \frac{\sigma_h}{N}, & \quad \forall K \in \mathcal{T}_h
\label{eq-1}
\\
\frac{1}{d} \text{tr}( (F_K')^T M_K F_K') = \text{det}( (F_K')^T M_K F_K')^{\frac 1 d},
&\quad \forall K \in \mathcal{T}_h
\label{ali-1}
\end{align}
where $d$ is the dimension of the domain $\Omega$, $M_K$ is the average of $M$ over $K$,
$F_K$ is the affine mapping from the reference element
$\hat{K}$ to element $K$, $F_K'$ denotes the Jacobian matrix of $F_K$, and
$\sigma_h = \sum_{K \in \mathcal{T}_h} |K| \text{det}(M_K)^{\frac 1 2}$.
Condition (\ref{eq-1}) requires the elements to have the same size in the metric $M$
while condition (\ref{ali-1}) requires that they be equilateral in the metric.
For a given metric tensor, various mesh strategies can be used to generate meshes approximately
satisfying (\ref{eq-1}) and (\ref{ali-1}), including the variational approach \cite{Hua01b,HR11},
Delaunay-type triangulation \cite{BGHLS97,BGM97,CHMP97,PVMZ97}, advancing front \cite{GS98},
bubble meshing \cite{YS00}, and combination of refinement, local modification, and
smoothing or node movement \cite{ABHFDV02,BH96,DVBFH02,HDBAFV00,RLSF04}.

\section{Numerical examples}
\label{SEC:numerics}

In this section we present five two-dimensional examples to verify the theoretical analysis
in the previous two sections. Since the non-obtuse angle condition (\ref{thm:irreducibleM-1})
is stronger than the Delaunay-type condition (\ref{thm:irreducibleM-2}) in 2D, we shall focus
on the latter in this section. For convenience, we define
\begin{align}
\alpha_{\max, \mathbb{D}^{-1}} & =  \max\limits_{K \in \mathcal{T}_h} \alpha_{\max, \mathbb{D}^{-1}}^K ,
\label{alpha-max}
\\
\alpha_{\text{sum}, \mathbb{D}^{-1}} & =  \max\limits_{e_{j,k} \in \mathcal{T}_h}
\frac{1}{2} \left [ \frac{}{} \alpha_{jk, \mathbb{D}^{-1}}^K + \alpha_{jk, \mathbb{D}^{-1}}^{K'}
+ \text{ arccot}\left (\sqrt{\frac{\mbox{ det}(\mathbb{D}_{K'})}{\mbox{ det}(\mathbb{D}_{K})}}
\cot (\alpha_{jk, \mathbb{D}^{-1}}^{K'} ) \right ) \right.
\nn \\
& \left .
 \qquad \qquad + \text{ arccot}\left (\sqrt{\frac{\mbox{ det}(\mathbb{D}_{K})}{\mbox{ det}(\mathbb{D}_{K'})}}
\cot (\alpha_{jk, \mathbb{D}^{-1}}^K ) \right ) \frac{}{} \right ]  .
\label{alpha-sum}
\end{align}
In our computation, principal eigenfunctions are normalized such that they have the maximum value one.
It is noted that analytical expressions for the principal eigenvalue and eigenfunctions are not available for
all of the examples. For convergence plot, we use a numerical principal eigenvalue obtained on a much finer
mesh as the reference value. We take $\Omega = (0,1)\times (0,1)$ in all but Example \ref{ex5.5} where
$\Omega = (0,1) \times (0,1)\backslash (\frac{4}{9}, \frac{5}{9})\times (\frac{4}{9}, \frac{5}{9})$.

\begin{figure}[thb]
\centering
\hspace{0.5cm}
\hbox{
\begin{minipage}[t]{2.5in}
\centerline{(a) Mesh135, $\alpha_{\text{sum}, \mathbb{D}^{-1}} = 1.71 \pi$}
\includegraphics[width=2in]{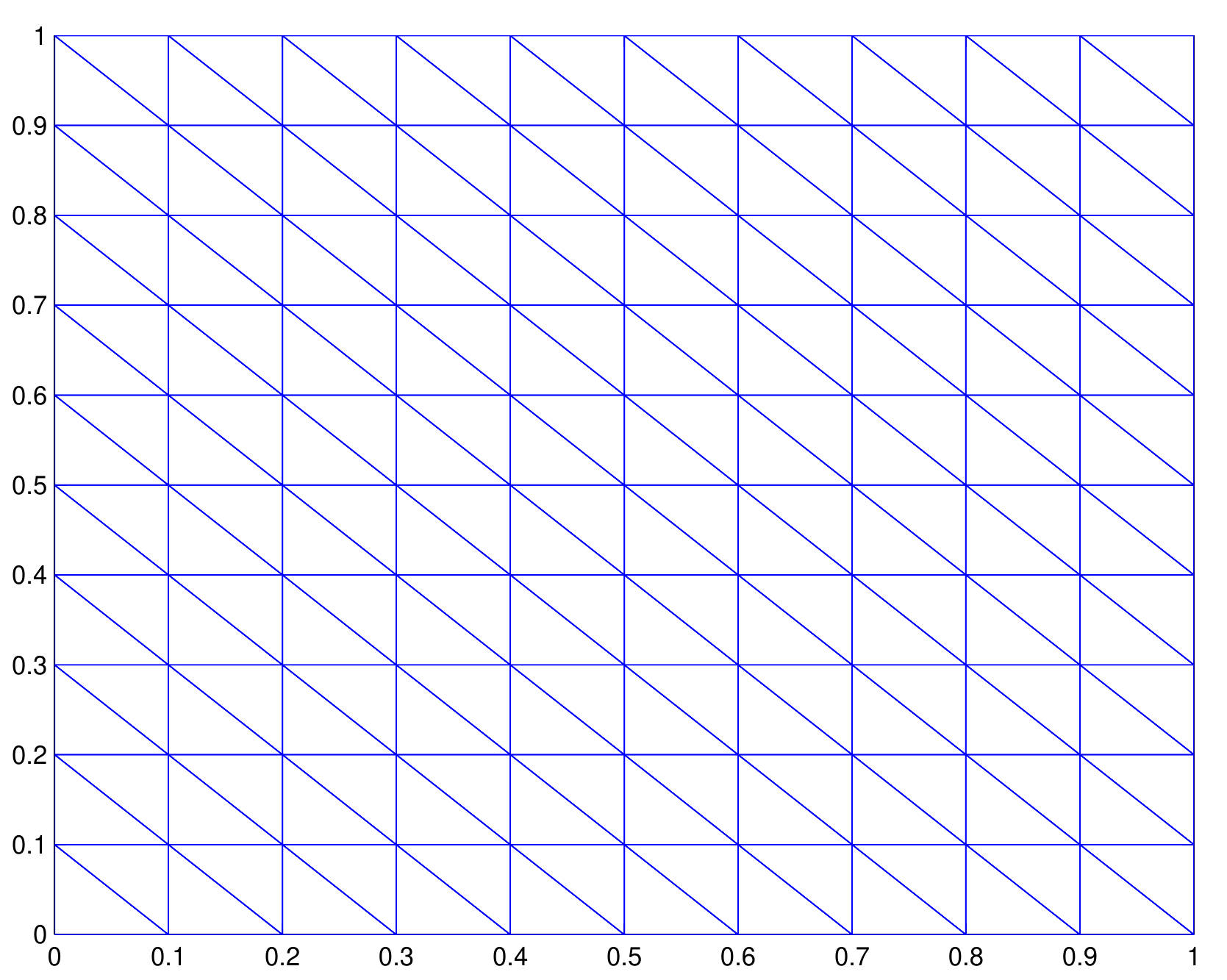}
\end{minipage}
\hspace{10mm}
\begin{minipage}[t]{2.5in}
\centerline{(b) Mesh45, $ \alpha_{\text{sum}, \mathbb{D}^{-1}} = 0.86 \pi$}
\includegraphics[width=2in]{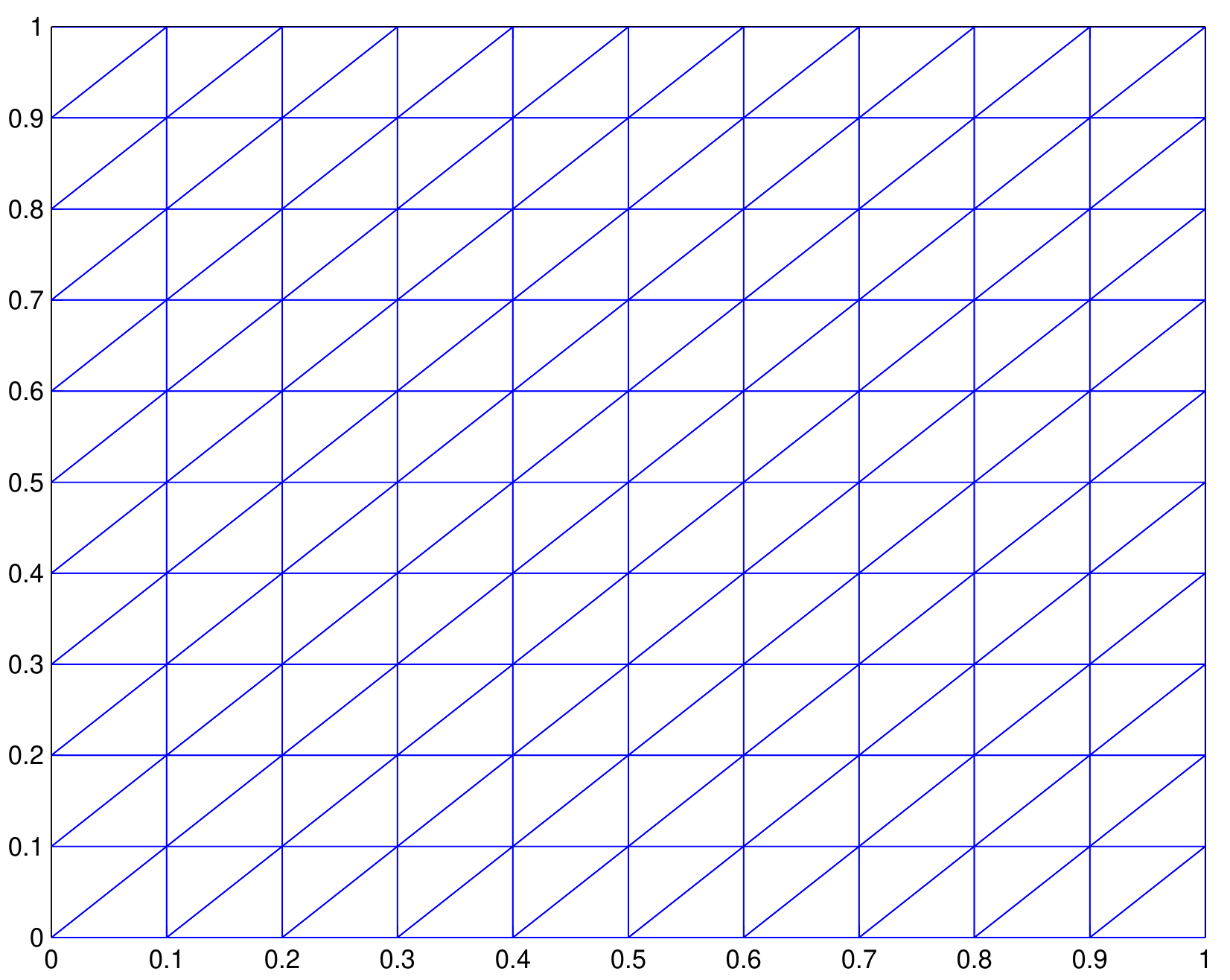}
\end{minipage}
}
\caption{Typical Mesh45 and Mesh135 meshes used in the computation for Example~\ref{ex5.1}.}
\label{ex5.1-f1}
\end{figure}

\begin{figure}[thb]
\centering
\hspace{0.5cm}
\hbox{
\begin{minipage}[t]{2.5in}
\centerline{(a) with Mesh135}
\includegraphics[width=2in]{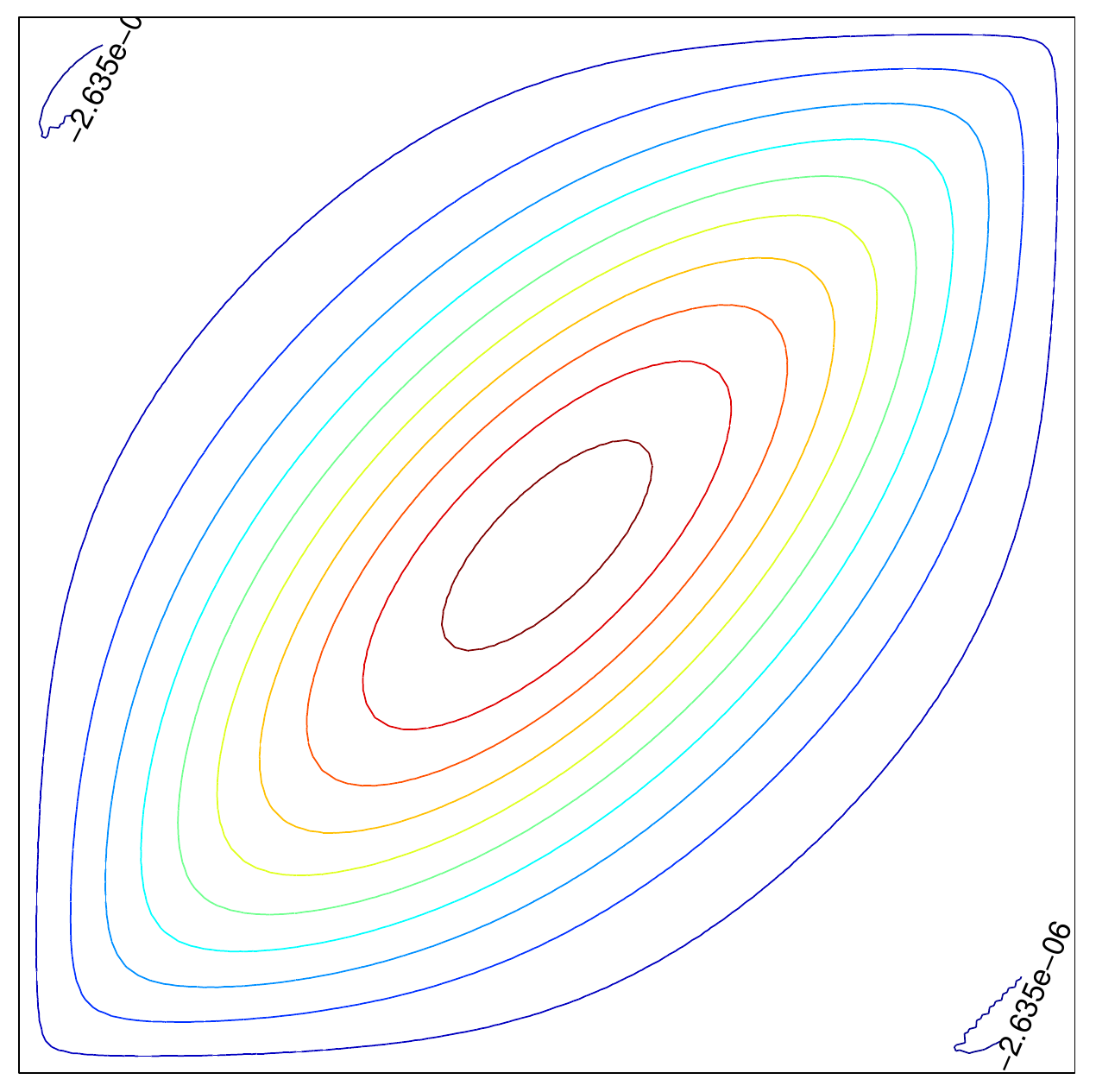}
\end{minipage}
\hspace{10mm}
\begin{minipage}[t]{2.5in}
\centerline{(b) with Mesh45}
\includegraphics[width=2in]{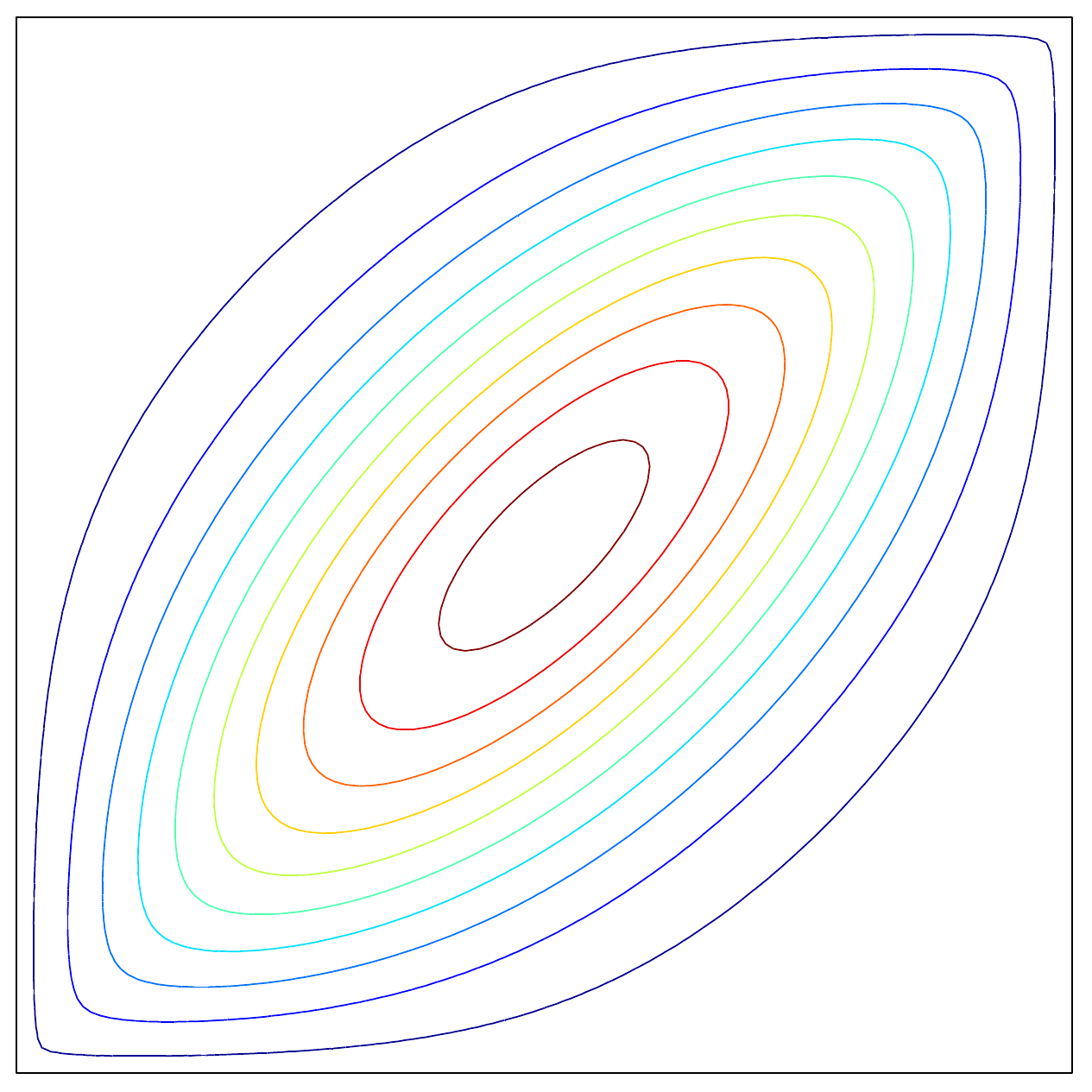}
\end{minipage}
}
\caption{Example~\ref{ex5.1}. Contours of the numerical eigenfunctions obtained with $J = 81$,
where $J$ is the number of mesh points in the $x$ (or $y$) axis.}
\label{ex5.1-f2}
\end{figure}

\begin{figure}[thb]
\centering
\hbox{
\begin{minipage}[t]{3in}
\centerline{(a) Undershoot with Mesh135}
\includegraphics[width=3in]{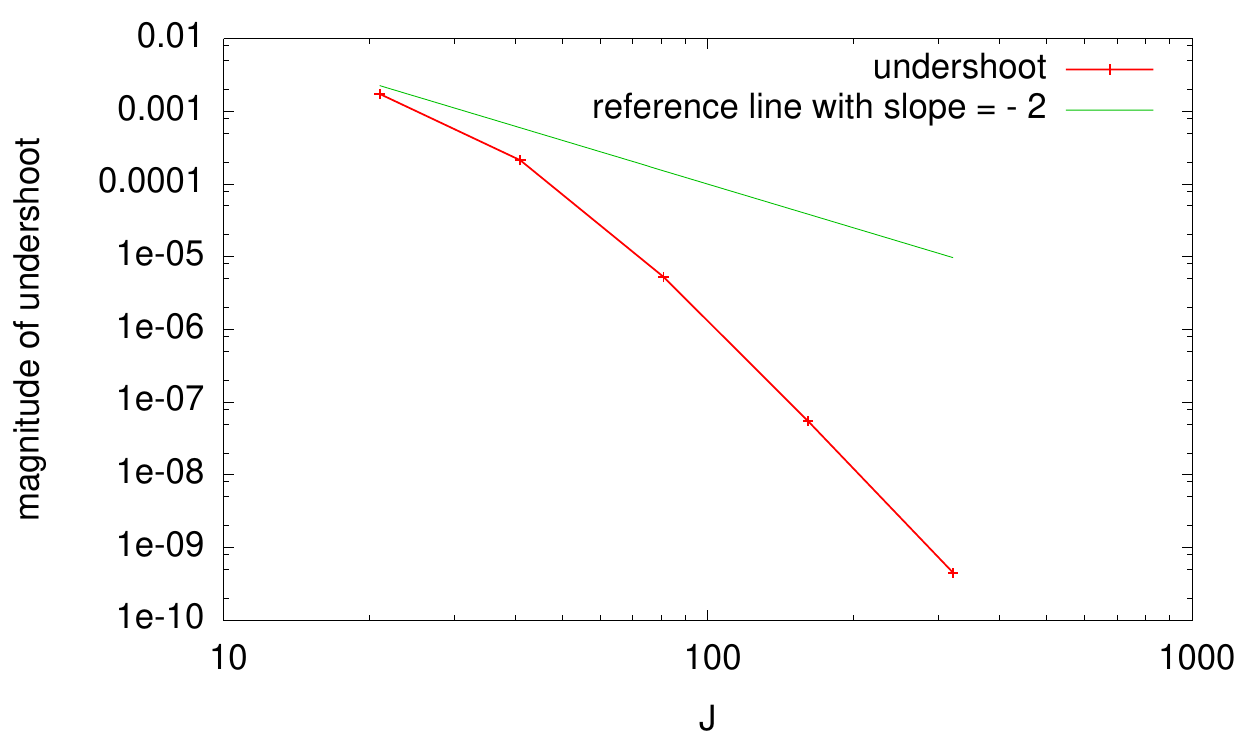}
\end{minipage}
\hspace{10mm}
\begin{minipage}[t]{3in}
\centerline{(b) Error in $\lambda_1$ with Mesh45 and Mesh135}
\includegraphics[width=3in]{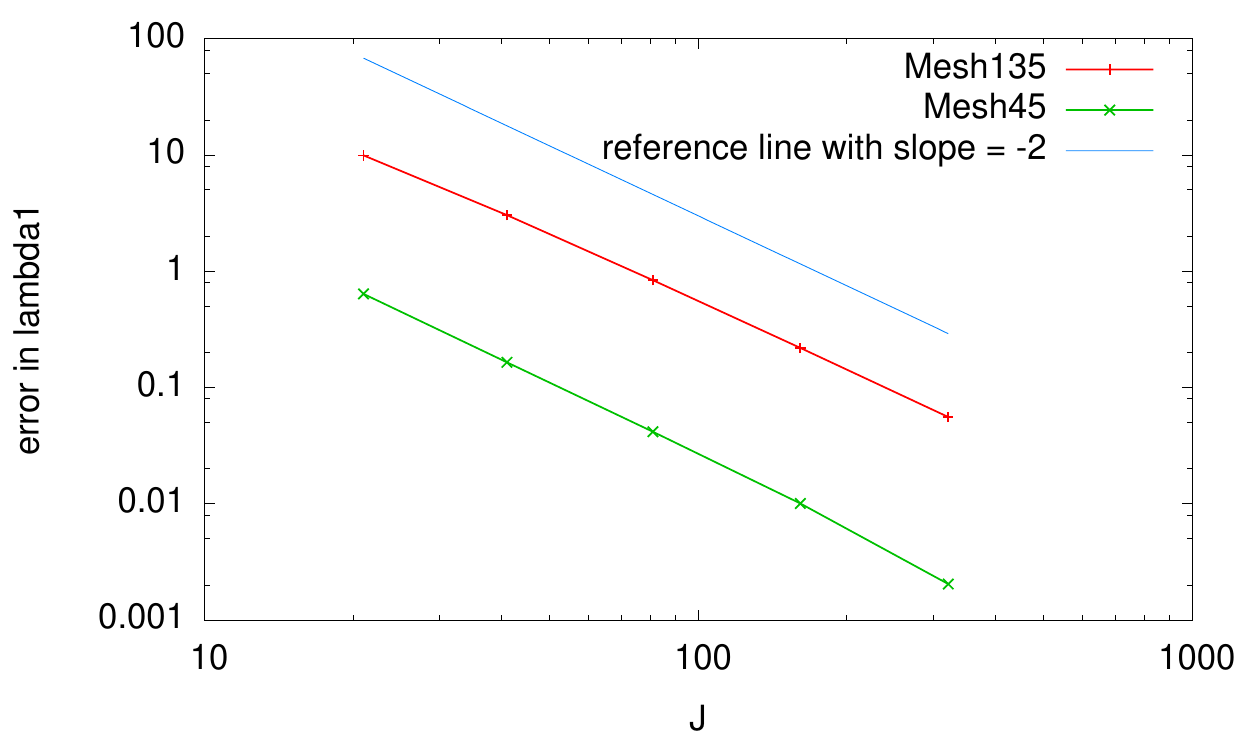}
\end{minipage}
}
\caption{Example~\ref{ex5.1}. The magnitude of the undershoot in the computed principal
eigenfunction and the error in the computed principal eigenvalue are plotted as function of $J$.
In (b), the reference value is $\lambda_1 \approx 150.288$ which is obtained with a Mesh45 of $J = 641$.}
\label{ex5.1-f3}
\end{figure}

\begin{exam}
\label{ex5.1}
The first example is in the form of (\ref{eigen-1}) with
\beq
\mathbb{D} = \left [\begin{array}{rr} 10 & 9 \\ 9 & 10 \end{array} \right ],\quad \V{b} = 0, \quad c = 0.
\label{ex5.1-1}
\eeq
Two types of mesh are used in the computation, Mesh135 and Mesh45.
They are obtained by cutting each square of a rectangular mesh into
two right triangles along the northwest or northeast diagonal line; see Fig~\ref{ex5.1-f1}.
For Mesh135, we have $\alpha_{\max, \mathbb{D}^{-1}} = 0.86 \pi$, 
$ \alpha_{\text{sum}, \mathbb{D}^{-1}} = 1.71 \pi$ and for Mesh45,
$\alpha_{\max, \mathbb{D}^{-1}} = 0.43 \pi$, $ \alpha_{\text{sum}, \mathbb{D}^{-1}} = 0.86 \pi$.
Thus, Mesh45 satisfies both (\ref{thm:irreducibleM-1}) and (\ref{thm:irreducibleM-2})
whereas Mesh135 does not satisfy any of them.

Fig.~\ref{ex5.1-f2} shows the contours of the numerical approximations of the principal eigenfunction
obtained with Mesh135 and Mesh45. It can be seen that the eigenfunction obtained with Mesh135
has some negative values (undershoot) near the southeast and northwest corners whereas the one with Mesh45
has no undershoot or overshoot. The magnitude of the undershoot is plotted in Fig.~\ref{ex5.1-f3}(a)
as the mesh is refined. The figure shows that the undershoot decreases at a rate much faster than
the approximation order (i.e., the second order for P1 linear finite elements) but never disappears
even for a fine mesh. Fig.~\ref{ex5.1-f3}(b) shows the second order convergence for the principal
eigenvalue for both types of mesh although the result with Mesh45 is a magnitude more accurate
than that with Mesh135.
\end{exam}

\begin{figure}[thb]
\centering
\hspace{0.5cm}
\hbox{
\begin{minipage}[t]{2.7in}
\centerline{(a) with Mesh135 of $J=41$}
\includegraphics[width=2.7in]{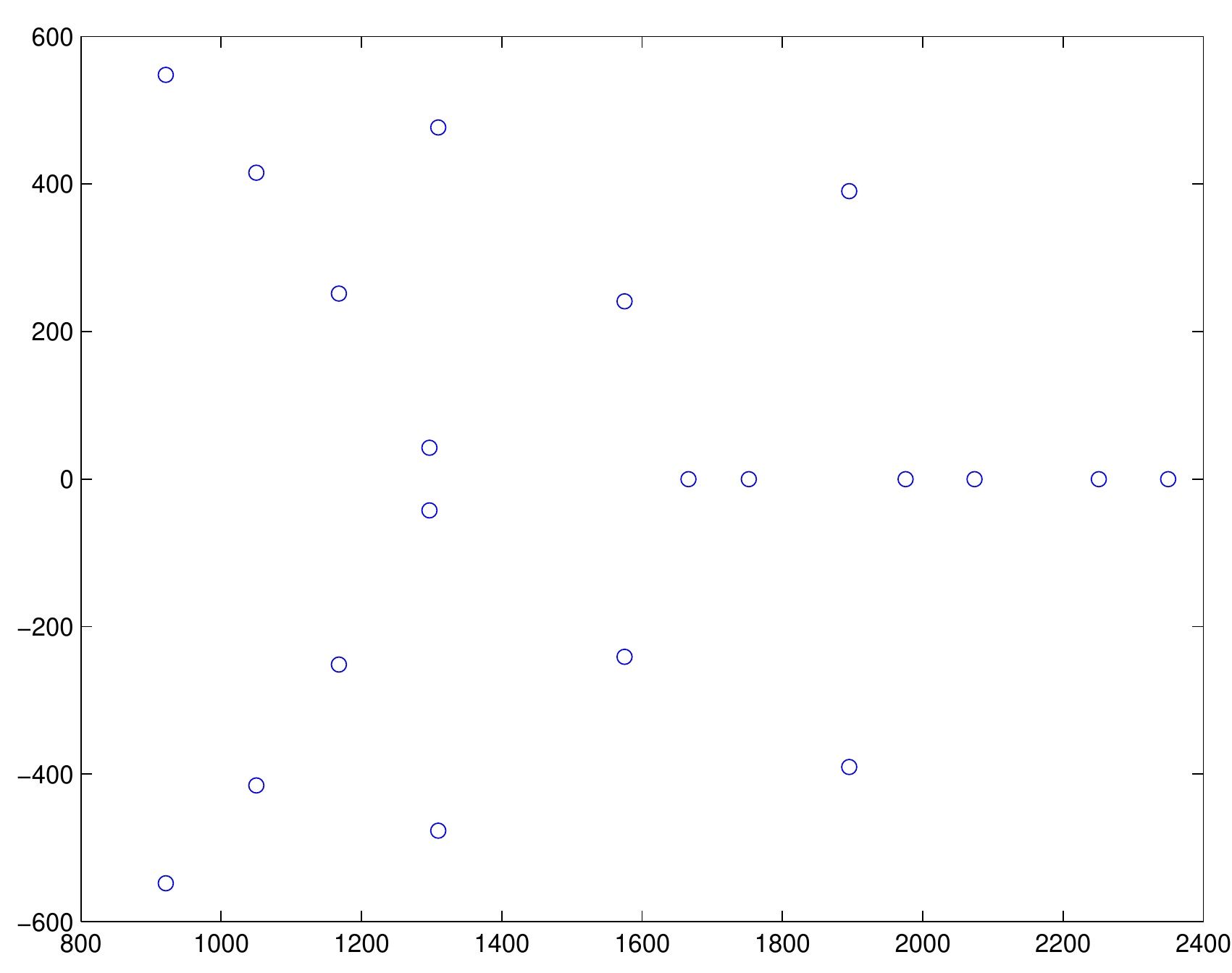}
\end{minipage}
\hspace{10mm}
\begin{minipage}[t]{2.7in}
\centerline{(b) with Mesh135 of $J=81$}
\includegraphics[width=2.7in]{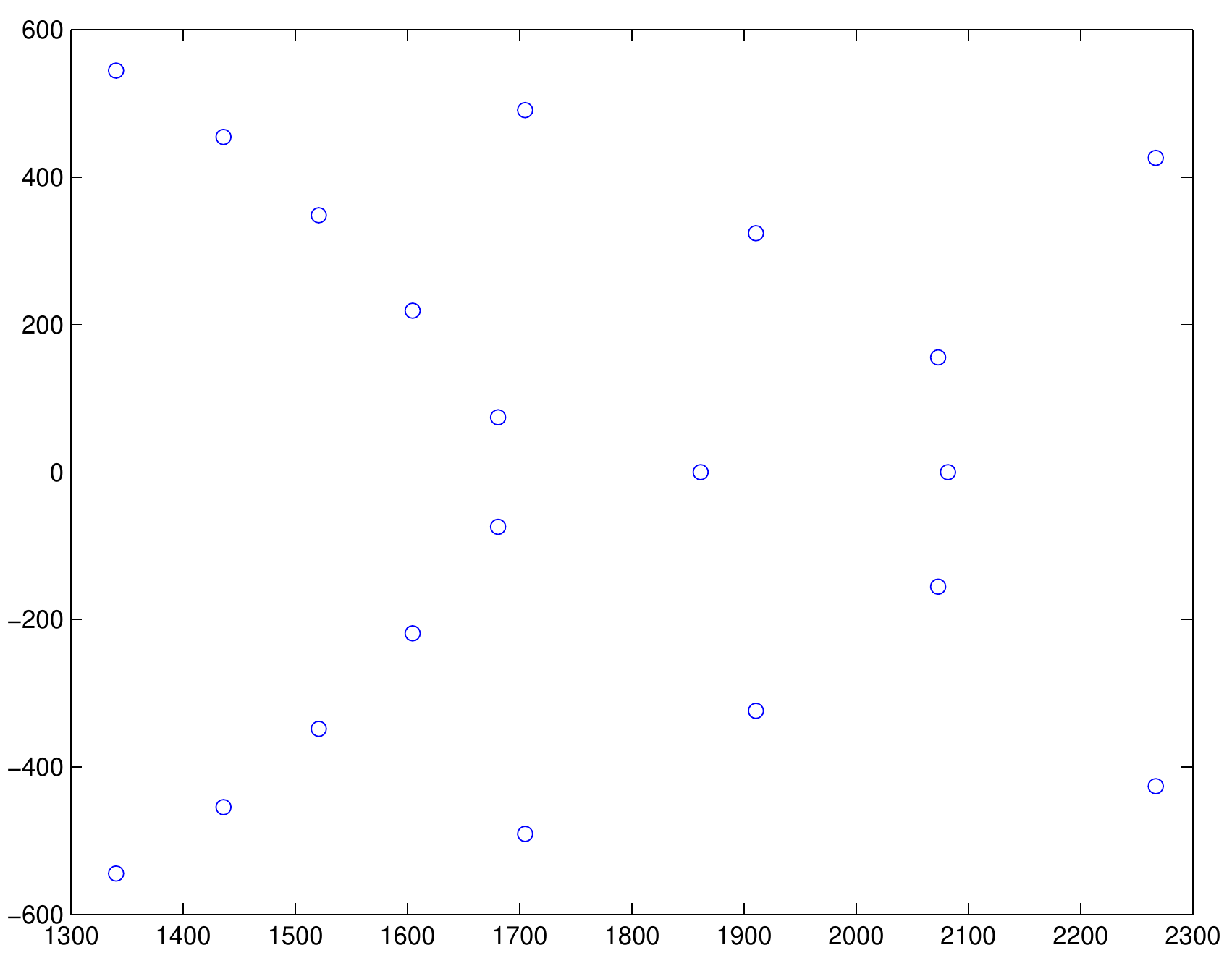}
\end{minipage}
}
\caption{Example~\ref{ex5.2}. The distribution of the first twenty smallest (in modulus) eigenvalues
for the discrete eigenvalue problem (\ref{fem-1}) with Mesh135.}
\label{ex5.2-f1}
\end{figure}

\begin{figure}[thb]
\centering
\hspace{0.5cm}
\hbox{
\begin{minipage}[t]{2.7in}
\centerline{(a) with Mesh45 of $J=41$}
\includegraphics[width=2.7in]{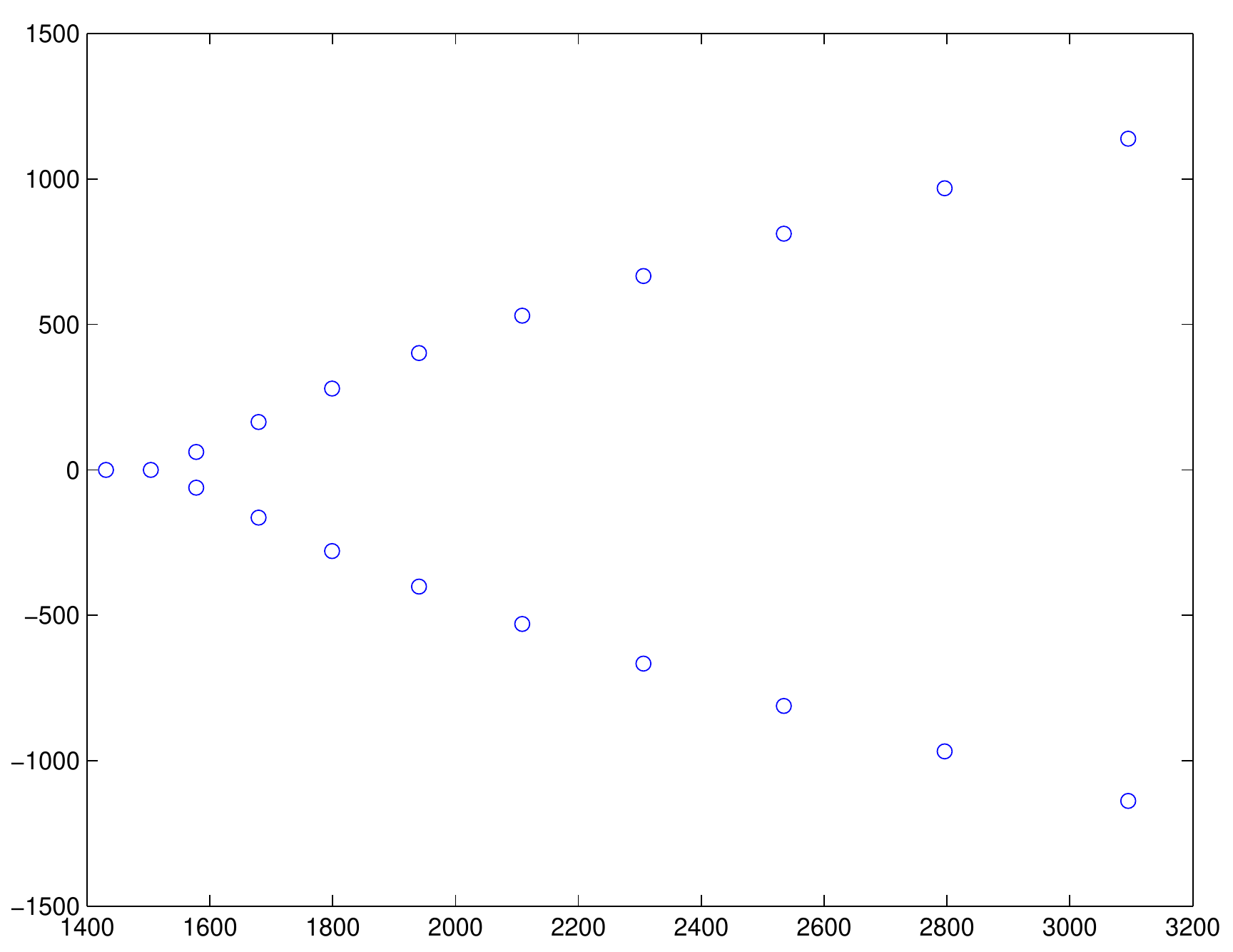}
\end{minipage}
\hspace{10mm}
\begin{minipage}[t]{2.7in}
\centerline{(b) with Mesh45 of $J=81$}
\includegraphics[width=2.7in]{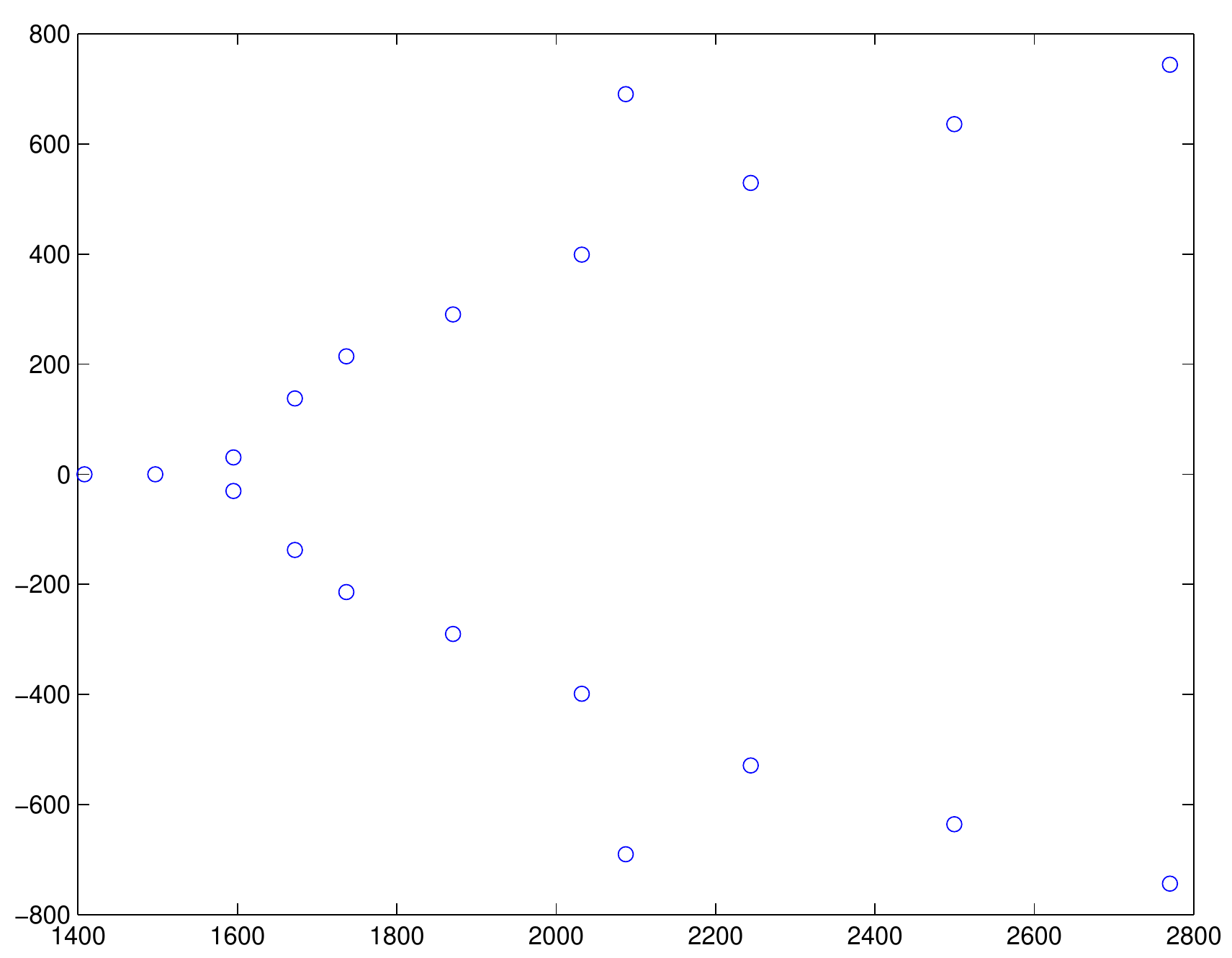}
\end{minipage}
}
\caption{Example~\ref{ex5.2}. The distribution of the first twenty smallest (in modulus) eigenvalues
for the discrete eigenvalue problem (\ref{fem-1}) with Mesh45.}
\label{ex5.2-f2}
\end{figure}

\begin{figure}[thb]
\centering
\hbox{
\begin{minipage}[t]{2.5in}
\centerline{(a) Principal eigenfunction}
\begin{center}
\includegraphics[width=1.6in]{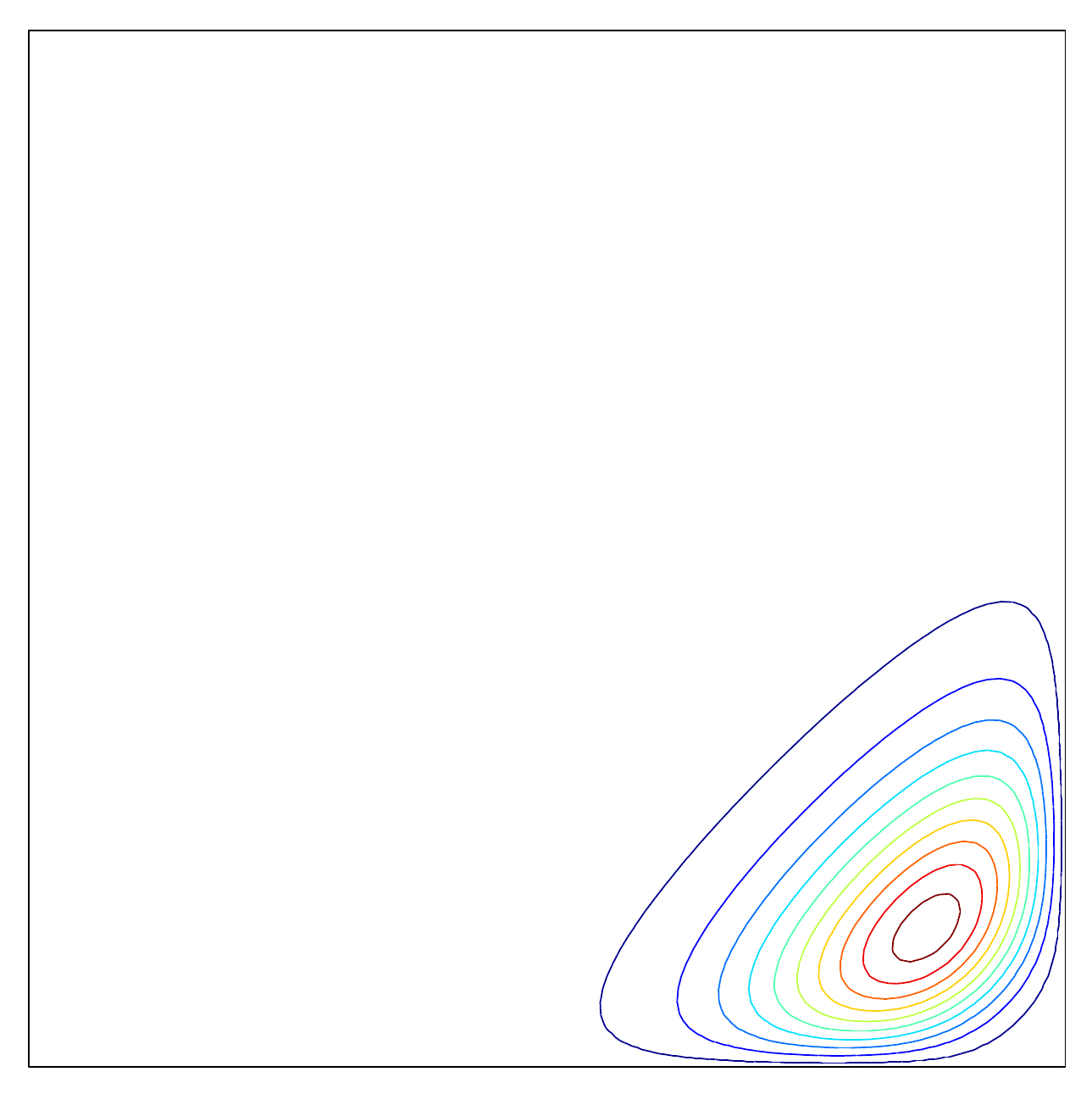}
\end{center}
\end{minipage}
\hspace{10mm}
\begin{minipage}[t]{3in}
\centerline{(b) Error in $\lambda_1$ (with the reference value $1401.39$)}
\includegraphics[width=3in]{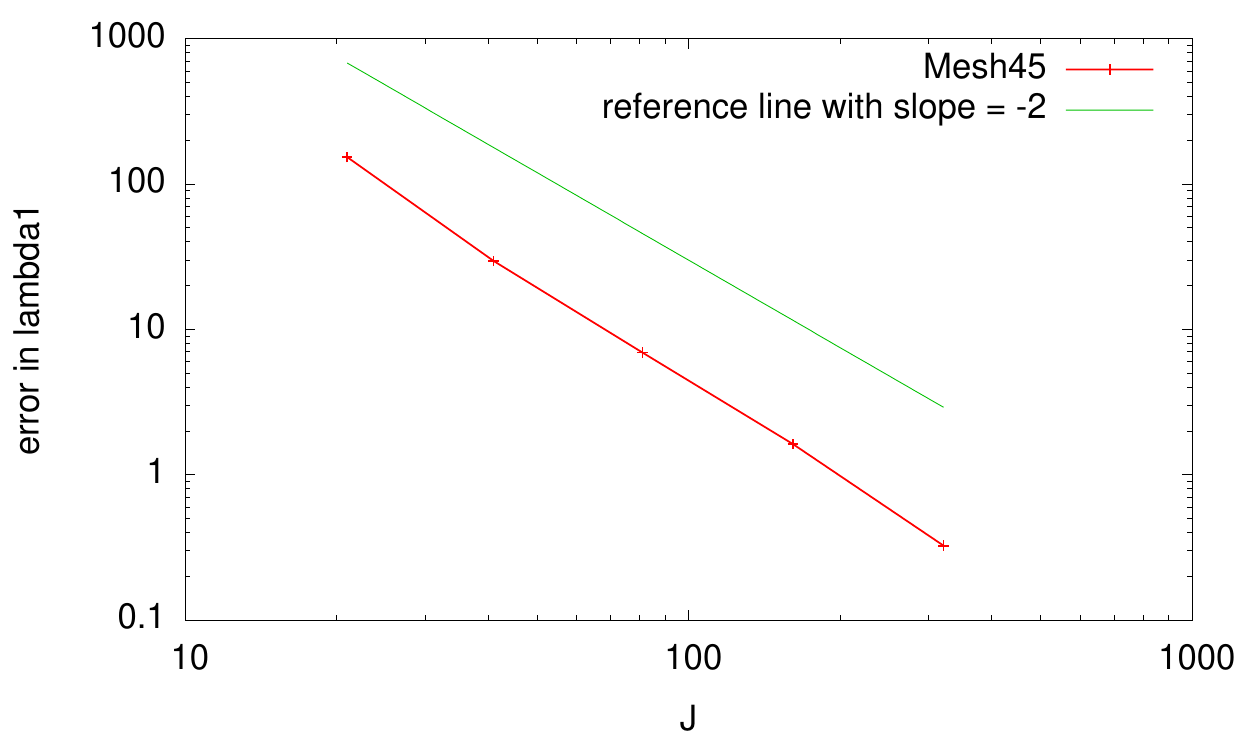}
\end{minipage}
}
\caption{Example~\ref{ex5.2}. (a) The contours of the numerical eigenfunction obtained with Mesh45
$J=81$ and (b) the error in the computed principal eigenvalue plotted as function of $J$.}
\label{ex5.2-f3}
\end{figure}

\begin{exam}
\label{ex5.2}
The second example is in the form of (\ref{eigen-1}) with
\beq
\mathbb{D} = \left [\begin{array}{rr} 10 & 9 \\ 9 & 10 \end{array} \right ],\quad \V{b} = \left [\begin{array}{c}50\\ -50\end{array}
\right ], \quad c = 1.
\label{ex5.2-1}
\eeq
Notice that this example is similar to the previous one except that this example contains both
the convection and reaction terms and is nonsymmetric. Both
Mesh135 and Mesh45 in Fig.~\ref{ex5.1-f1} are used in the computation.

Recall that Mesh135 does not satisfy the mesh condition (\ref{thm:irreducibleM-2}). The distribution of the
first twenty smallest (in modulus) eigenvalues obtained with Mesh135 ($J= 41$ and $J=81$)
is shown in Fig.~\ref{ex5.2-f1}. One can see that the smallest eigenvalues are actually complex.
On the other hand, Mesh45 satisfies (\ref{thm:irreducibleM-2}) when it is sufficiently fine.
The distribution of the first twenty smallest eigenvalues obtained with Mesh45 ($J= 41$ and $81$)
is shown in Fig.~\ref{ex5.2-f2}. It can be seen that the smallest eigenvalue of the discrete problem
(\ref{fem-1}) is real for both cases with $J = 41$ and $81$. Moreover, the figure shows that
$\text{Re}(\lambda^h)\ge \lambda_1$ at least for the first twenty smallest eigenvalues.

Fig.~\ref{ex5.2-f3}(a) shows the contours of
the numerical eigenfunction obtained with Mesh45 $J=81$ and
Fig.~\ref{ex5.2-f3}(b) shows the error in the computed principal eigenvalue as function of $J$.
\end{exam}

\begin{figure}[thb]
\centering
\includegraphics[width=2.5in]{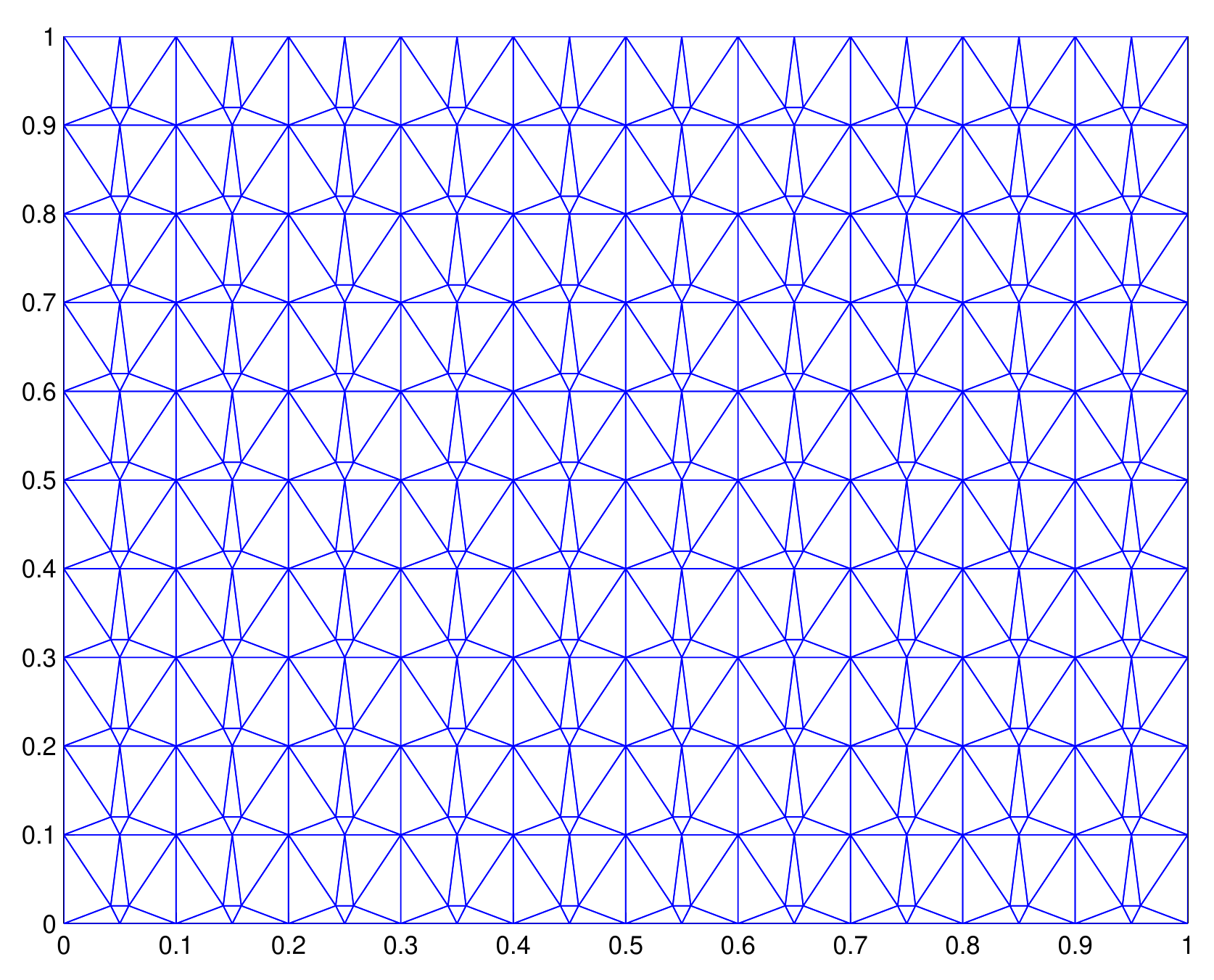}
\caption{An example ($J = 11$) of the mesh  used in the computation for Example~\ref{ex5.3} is shown, with
$\alpha_{\max, \mathbb{D}^{-1}} = 0.49 \pi$.}
\label{ex5.3-f1}
\end{figure}

\begin{figure}[thb]
\centering
\hbox{
\begin{minipage}[t]{2.5in}
\centerline{(a) Contours of computed eigenfunction}
\begin{center}
\includegraphics[width=1.6in]{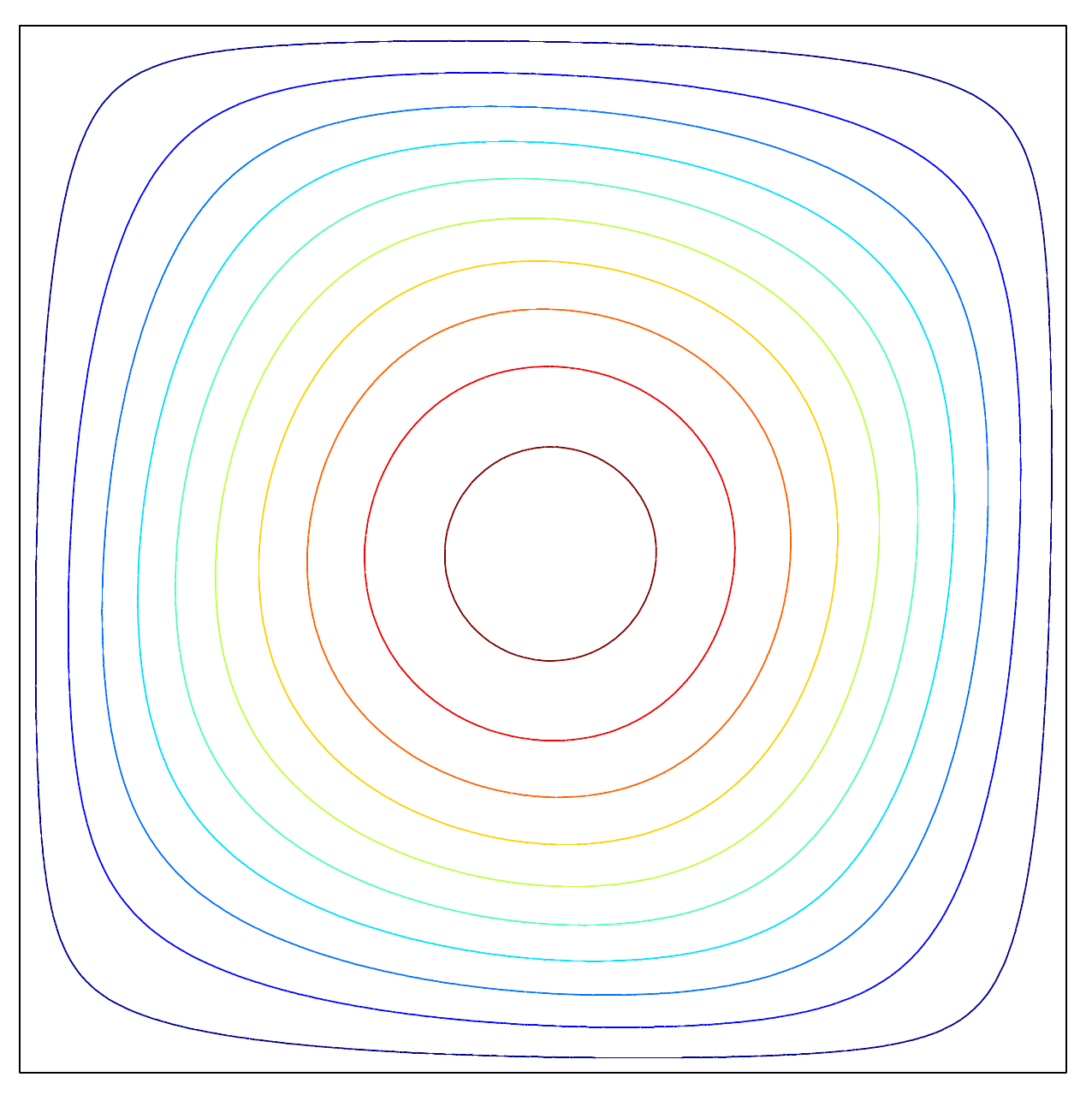}
\end{center}
\end{minipage}
\hspace{10mm}
\begin{minipage}[t]{3in}
\centerline{(b) Error in $\lambda_1$ (with the reference value $21.0714$)}
\begin{center}
\includegraphics[width=3in]{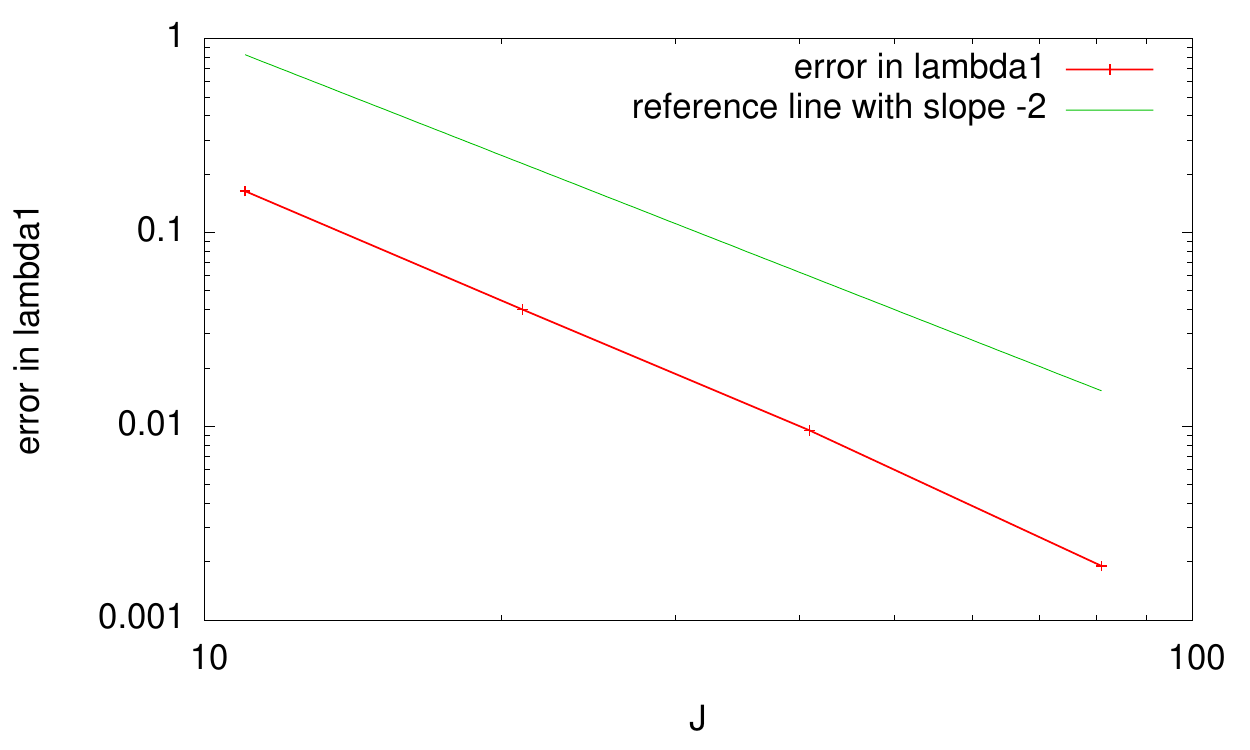}
\end{center}
\end{minipage}
}
\caption{Example~\ref{ex5.3}. (a) The contours of the computed principal eigenfunction with $J = 81$
and (b) the error in $\lambda_1$ plotted as function of $J$.}
\label{ex5.3-f2}
\end{figure}

\begin{exam}
\label{ex5.3}
The next example is in the form of (\ref{eigen-1}) with
\beq
\mathbb{D} = I + 0.05 \left [\begin{array}{cc} \cos(\pi x) & 0 \\ 0 & \sin(\pi y) \end{array} \right ],\quad
\V{b} = \left [\begin{array}{c} 20 (y-0.5) \\ -20 (x-0.5)\end{array} \right ], \quad c = 1.
\label{ex5.3-1}
\eeq
Notice that the diffusion matrix and the convection vector are functions of $x$ and $y$. The diffusion
matrix is chosen as a small perturbation of the identity matrix so that the acute mesh (in the Euclidean
sense) shown in Fig.~\ref{ex5.3-f1} is also acute in the metric specified by $\mathbb{D}^{-1}$.
As a consequence, the mesh condition (\ref{thm:irreducibleM-1}) (and therefore (\ref{thm:irreducibleM-2}))
can be satisfied when the mesh is sufficiently fine.

The contours of a computed principal eigenfunction (with $J = 81$) is shown in Fig.~\ref{ex5.3-f2}(a).
No undershoot is observed in the solution. The error in $\lambda_1$ is plotted in Fig.~\ref{ex5.3-f2}(b)
as a function of $J$. The convergence rate is second order.
\end{exam}

\begin{figure}[thb]
\centering
\hbox{
\begin{minipage}[t]{2.5in}
\centerline{(a) Mesh, $N = 1129$, $\alpha_{\text{sum}, \mathbb{D}^{-1}} = 0.99 \pi$}
\begin{center}
\includegraphics[width=2in]{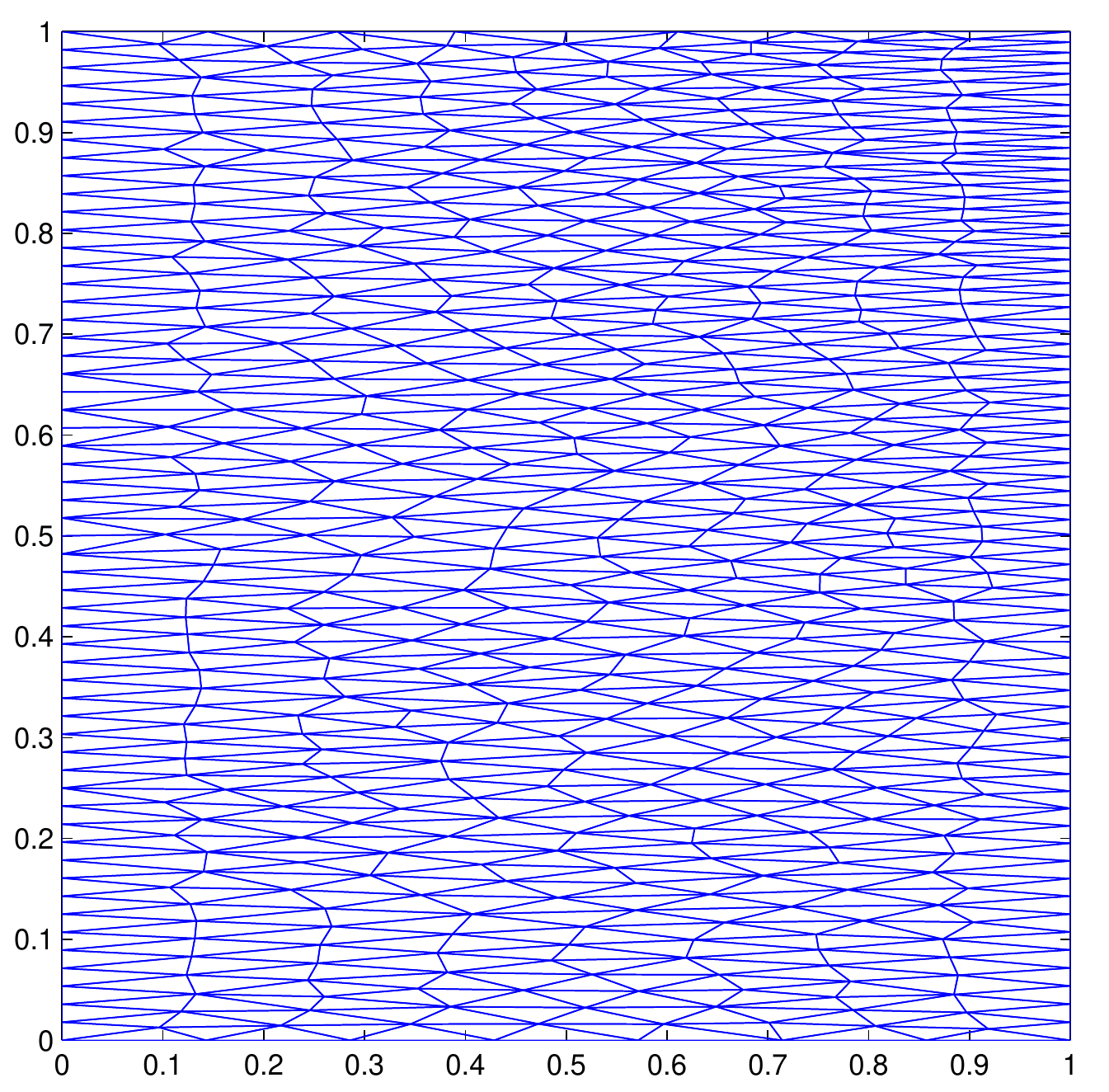}
\end{center}
\end{minipage}
\hspace{10mm}
\begin{minipage}[t]{3in}
\centerline{(b) Error in $\lambda_1$ (with the reference value $687.666$)}
\begin{center}
\includegraphics[width=3in]{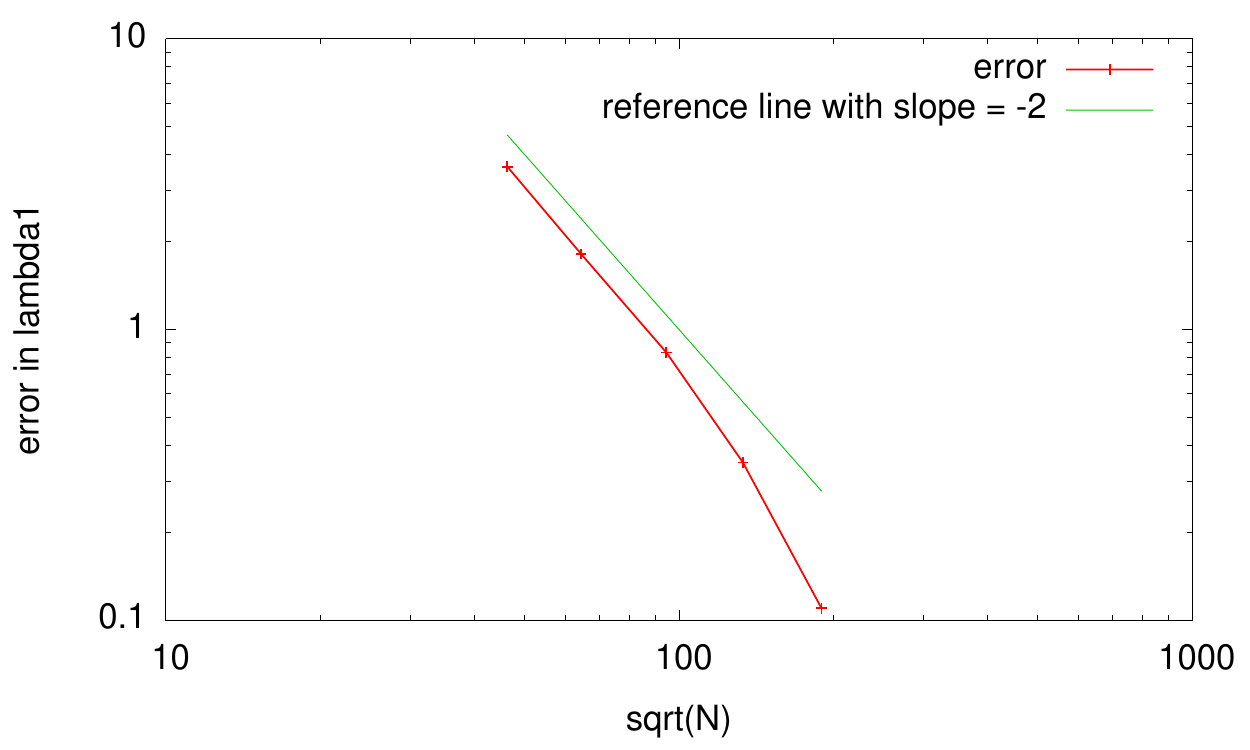}
\end{center}
\end{minipage}
}
\caption{Example~\ref{ex5.4}. A typical mesh used in the computation
and the convergence history in approximating
$\lambda_1$.}
\label{ex5.4-f1}
\end{figure}

\begin{figure}[thb]
\centering
\hspace{0.5cm}
\hbox{
\begin{minipage}[t]{2.5in}
\centerline{(a) contour plot}
\begin{center}
\includegraphics[width=2in]{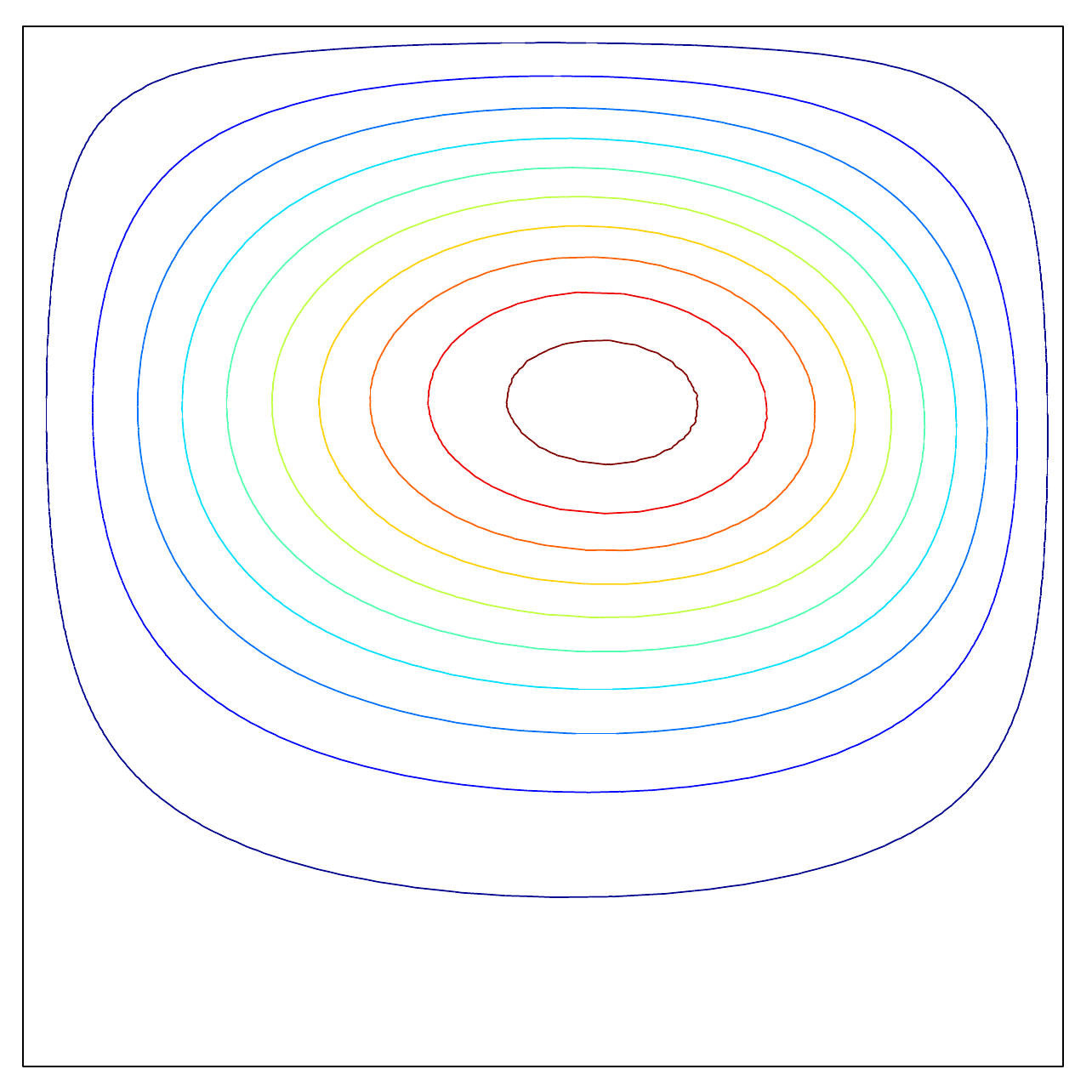}
\end{center}
\end{minipage}
\hspace{10mm}
\begin{minipage}[t]{2.5in}
\centerline{(b) surface plot}
\begin{center}
\includegraphics[width=2.5in]{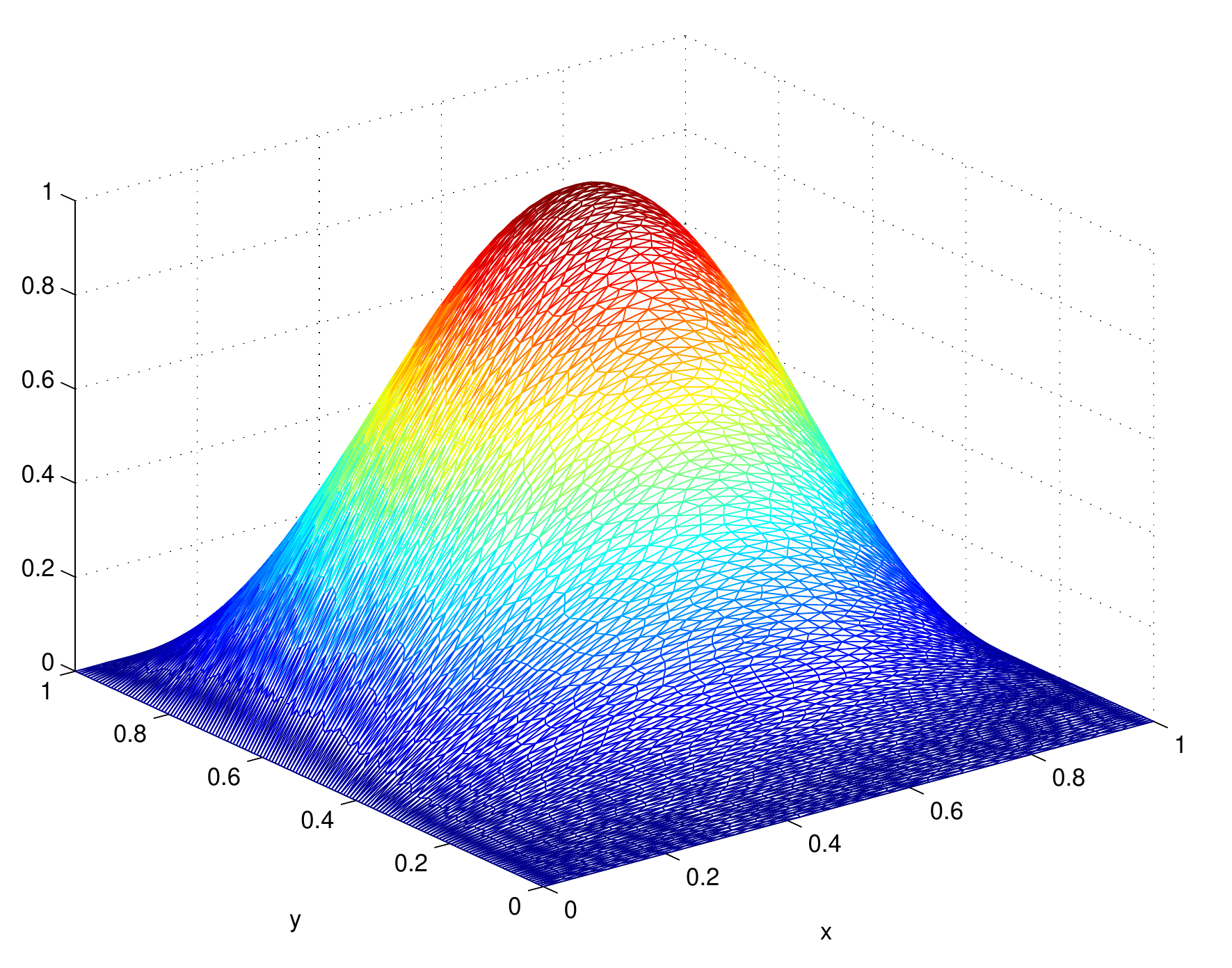}
\end{center}
\end{minipage}
}
\caption{Example~\ref{ex5.4}. Contour and surface plots for a computed principal eigenfunction obtained with $N=8925$.}
\label{ex5.4-f2}
\end{figure}

\begin{exam}
\label{ex5.4}
We have so far considered examples with constant or almost constant diffusion matrices.
In this and next examples, we consider the situation with variable diagonal and full diffusion matrices,
respectively. 

This example is in the form of (\ref{eigen-1}) with
\beq
\mathbb{D} = \left [ \begin{array}{cc} 100 (1-0.5\sin(x y \pi)) & 0 \\ 0 & (1+0.5 \cos(x y \pi)) \end{array}\right ],
\quad \V{b} = 0, \quad c = 0.
\label{ex5.4-1}
\eeq
Since $\mathbb{D}$ changes with location, it is impossible in general to predefine a mesh satisfying
the mesh condition (\ref{thm:irreducibleM-1}) or (\ref{thm:irreducibleM-2}).
We use here the BAMG (bidimensional anisotropic mesh generator) code developed
by Hecht \cite{Hec97} to generate approximate $M$-uniform meshes
for the metric tensor $M = \mathbb{D}^{-1}$ (cf. the discussion right after Theorem~\ref{thm:irreducibleM}).
BAMG is a Delaunay-type mesh generator \cite{CHMP97}
and allows the user to supply a metric tensor defined on a background mesh. It is used in our computation
in an iterative fashion: Starting from a coarse mesh, the metric tensor $M = \mathbb{D}^{-1}$
is computed and used in BAMG to generate a new mesh. The process is repeated ten times.

It is noted that $\mathbb{D}$ defined in (\ref{ex5.4-1}) is diagonal but very anisotropic, with the maximum
ratio of the two eigenvalues being over 100. Numerical results show that BAMG is able to generate
meshes satisfying (\ref{thm:irreducibleM-2}). Fig. \ref{ex5.4-f1}(a) shows such a mesh with
$\alpha_{\text{sum}, \mathbb{D}^{-1}} = 0.99 \pi$. No undershoot is observed in the computed
principal eigenfunctions, as shown in Fig.~\ref{ex5.4-f2}. A second order convergence rate
in approximating $\lambda_1$ is observed in Fig. \ref{ex5.4-f1}(b).
\end{exam}

\begin{figure}[thb]
\centering
\hspace{0.5cm}
\hbox{
\begin{minipage}[t]{2.5in}
\centerline{(a) Mesh with $N = 2260$, $\alpha_{\text{sum}, \mathbb{D}^{-1}} = 0.99 \pi$}
\begin{center}
\includegraphics[width=2in]{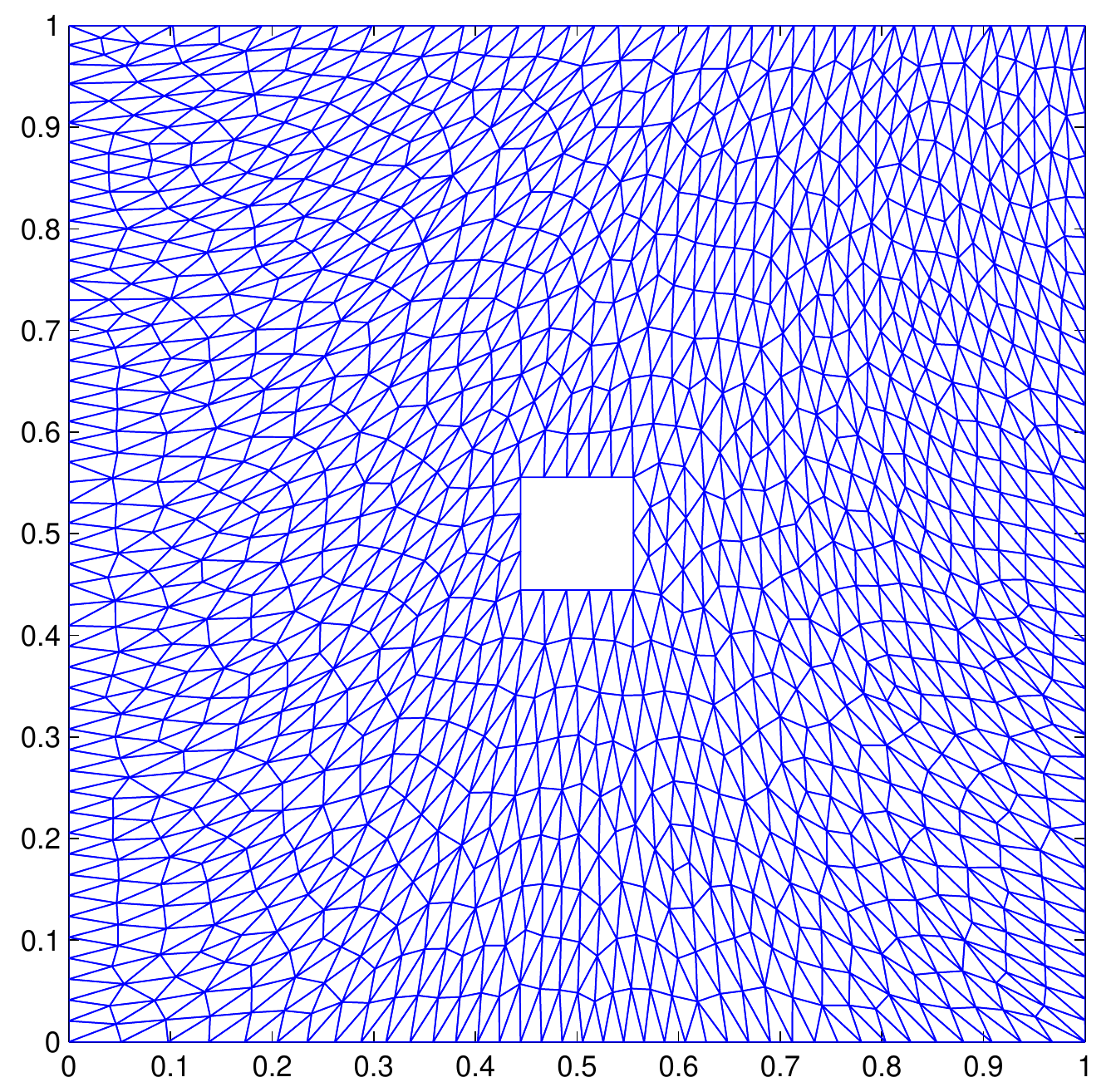}
\end{center}
\end{minipage}
\hspace{10mm}
\begin{minipage}[t]{2.5in}
\centerline{(b) Principal eigenfunction}
\begin{center}
\includegraphics[width=2in]{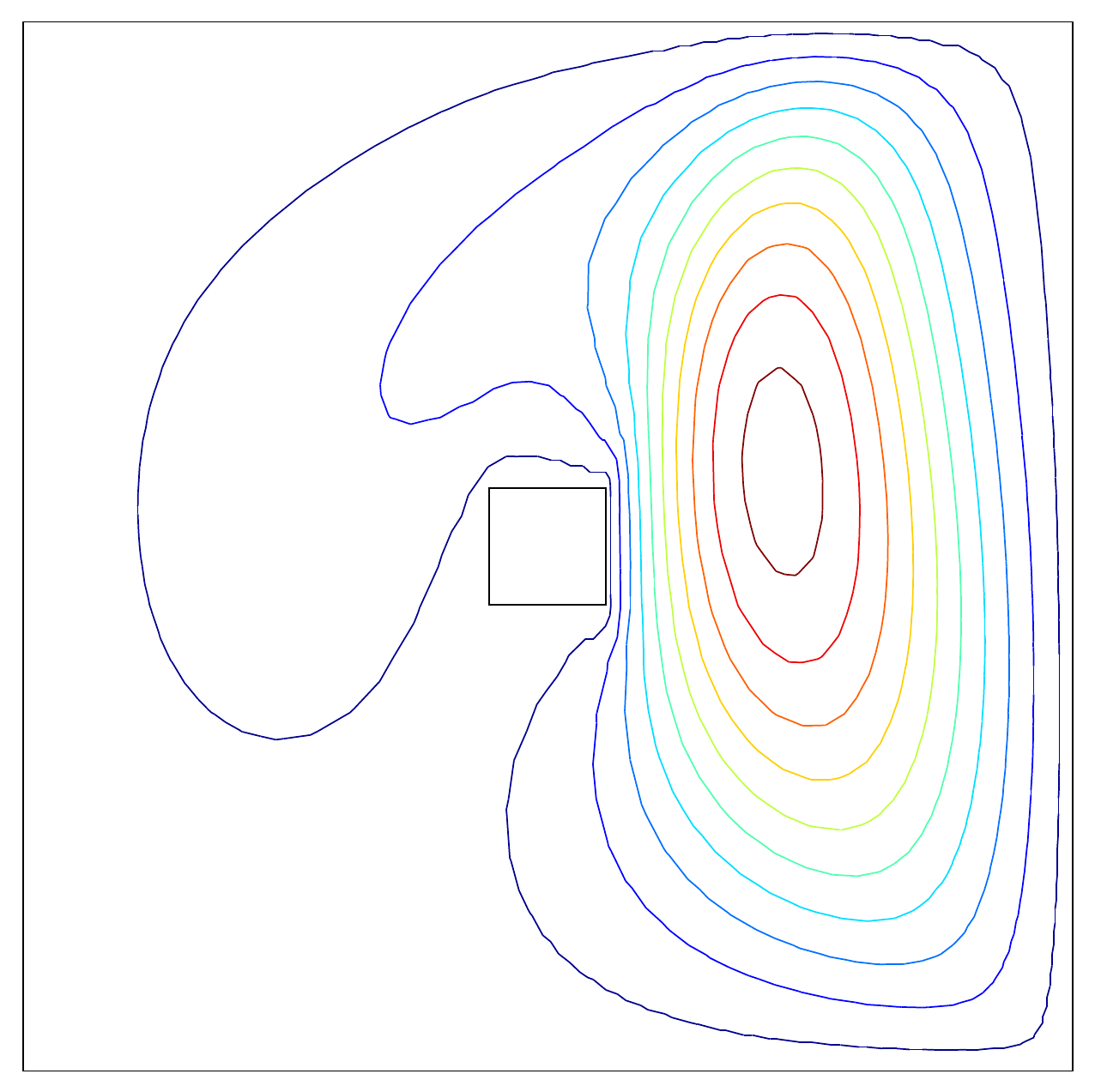}
\end{center}
\end{minipage}
}
\caption{Example~\ref{ex5.5}. A typical mesh generated with BAMG for metric tensor $M = \mathbb{D}^{-1}$ ($k=10$)
and the corresponding computed principal eigenfunction.}
\label{ex5.5-f1}
\end{figure}

\begin{figure}[thb]
\centering
\hspace{0.5cm}
\hbox{
\begin{minipage}[t]{2.5in}
\centerline{(a) Mesh with $N = 2373$, $\alpha_{\text{sum}, \mathbb{D}^{-1}} = 1.23 \pi$}
\begin{center}
\includegraphics[width=2in]{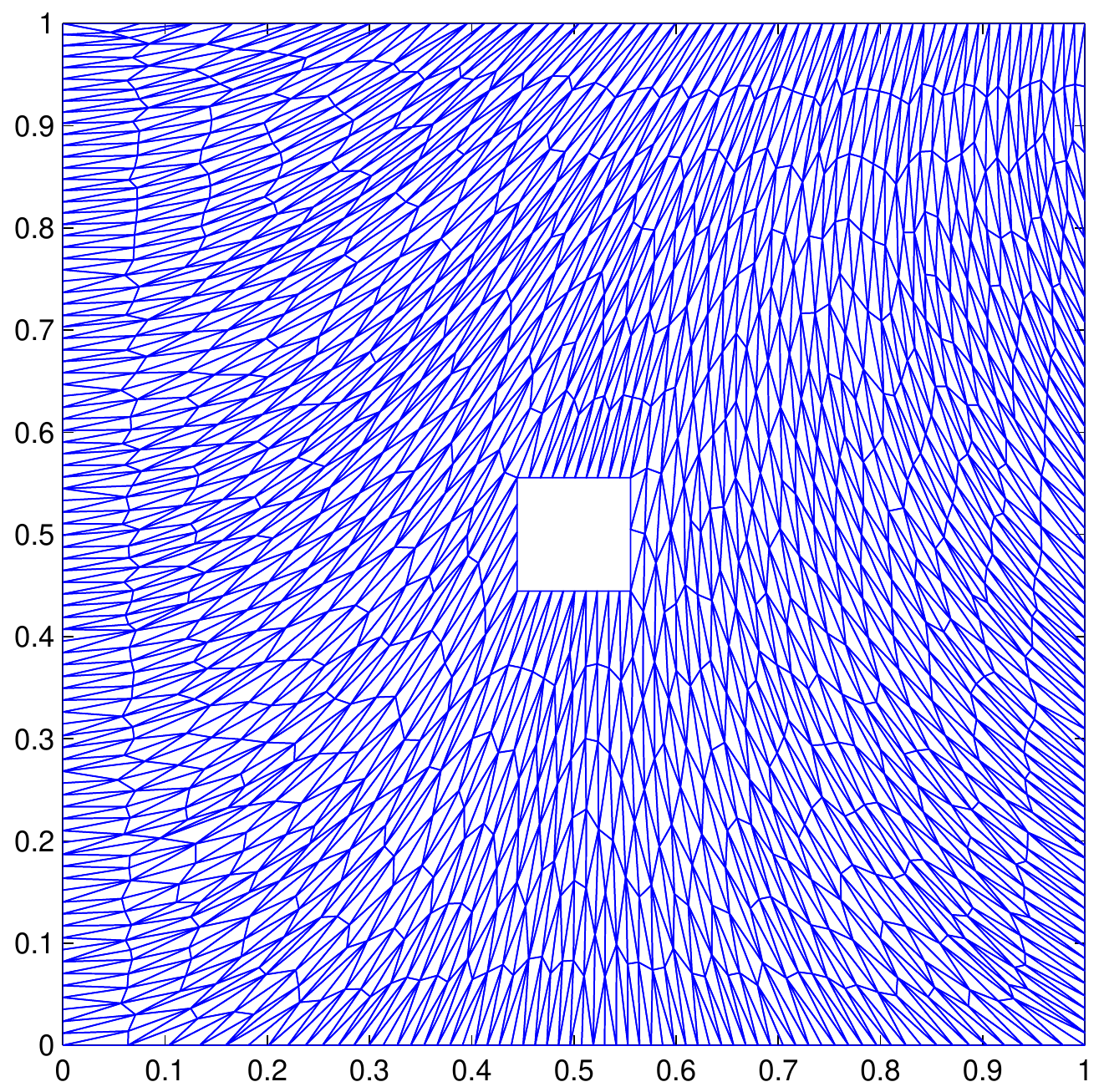}
\end{center}
\end{minipage}
\hspace{10mm}
\begin{minipage}[t]{2.5in}
\centerline{(b) Principal eigenfunction}
\begin{center}
\includegraphics[width=2in]{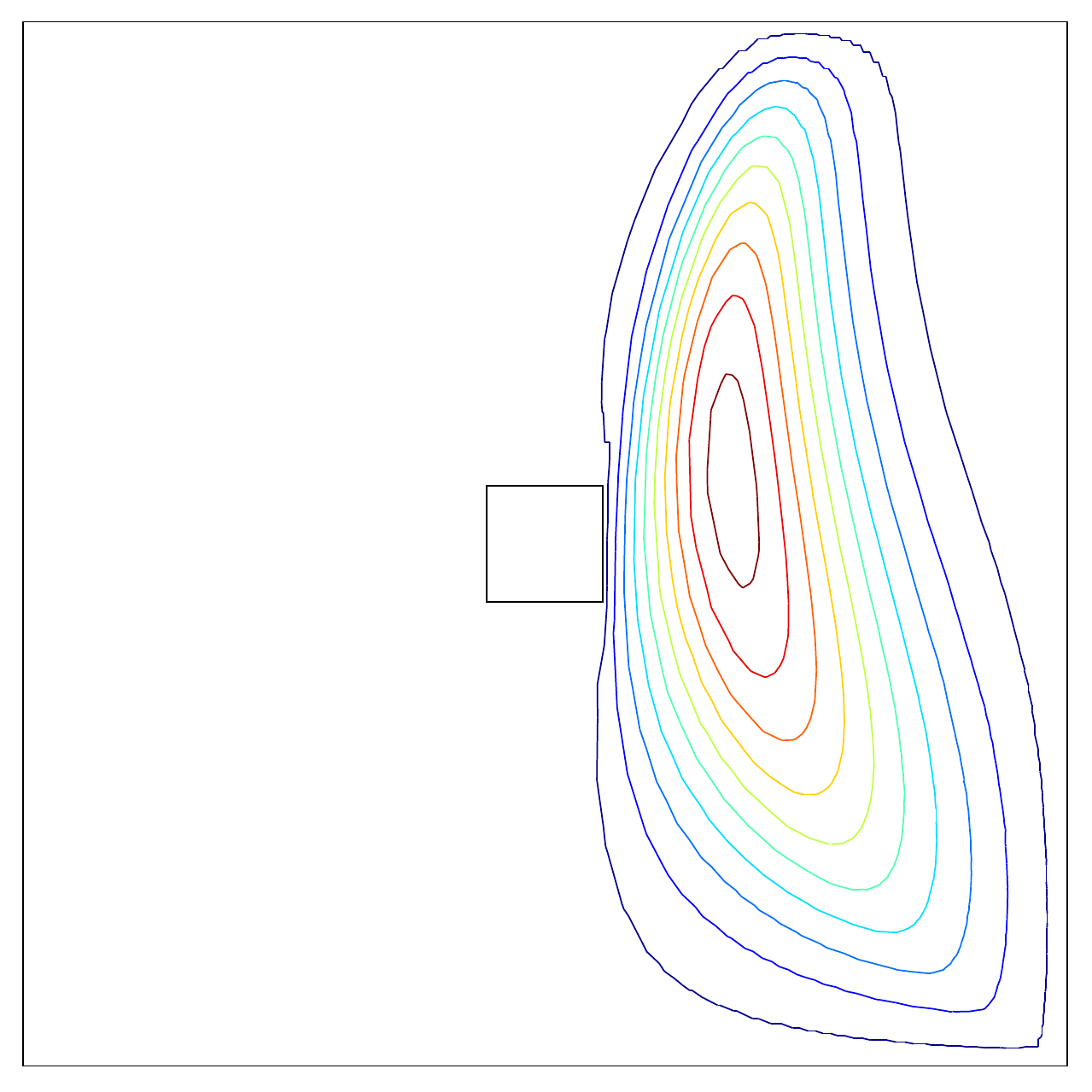}
\end{center}
\end{minipage}
}
\caption{Example~\ref{ex5.5}. A typical mesh generated with BAMG for metric tensor $M = \mathbb{D}^{-1}$ ($k=100$)
and the corresponding computed principal eigenfunction.}
\label{ex5.5-f2}
\end{figure}

\begin{figure}[thb]
\centering
\includegraphics[width=3in]{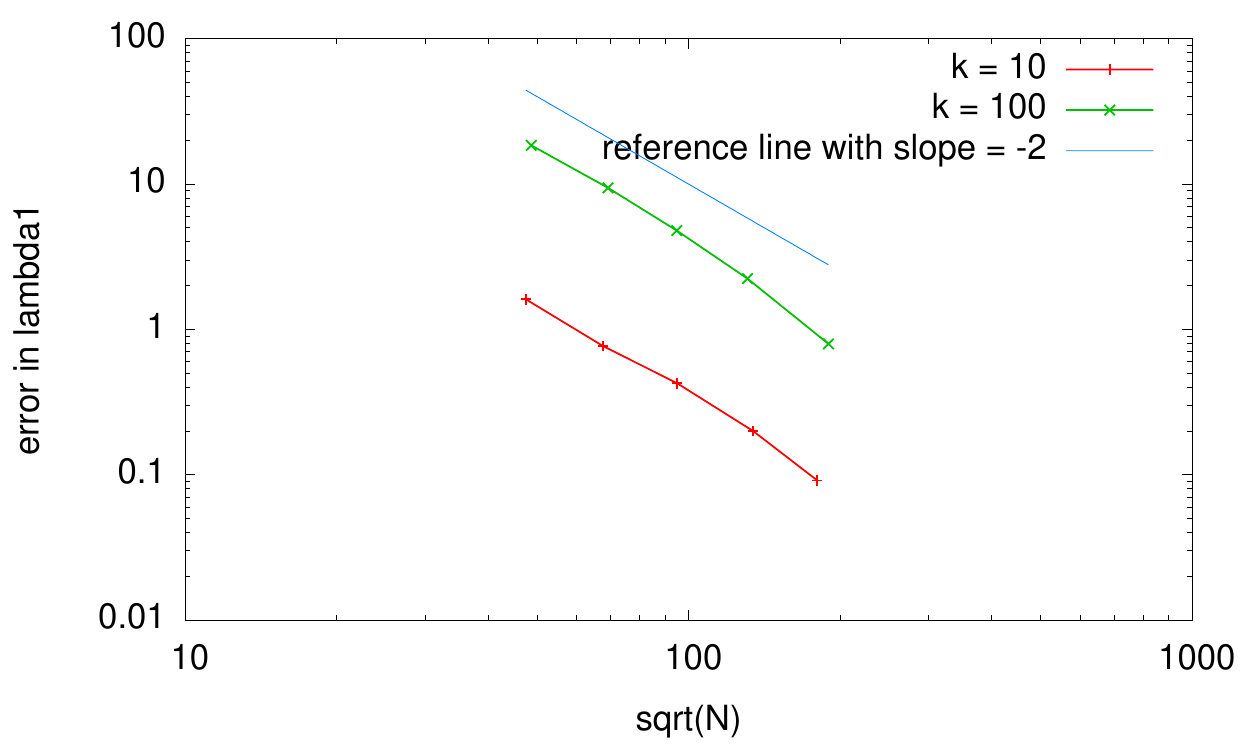}
\caption{Example~\ref{ex5.5}. The error in $\lambda_1$ is plotted as function of $\sqrt{N}$, where $N$ is the number
of mesh elements. The reference values for $\lambda_1$ are 170.422 for $k=10$ and 1020.15 for $k=100$.}
\label{ex5.5-f3}
\end{figure}

\begin{exam}
\label{ex5.5}
In this final example we consider a full diffusion matrix,
\begin{align}
& \mathbb{D} =  \left [\begin{array}{cc} \cos(\theta) & -\sin(\theta) \\ \sin(\theta) & \cos(\theta) \end{array} \right ]
\cdot 
\left [\begin{array}{cc} k (1-0.5\sin(x)\sin(y)) & 0 \\ 0 & (1+0.5\cos(x)\cos(y)) \end{array} \right ]
\nn \\
& \qquad \qquad \qquad  \cdot \left [\begin{array}{cc} \cos(\theta) & \sin(\theta) \\ -\sin(\theta) & \cos(\theta) \end{array} \right ],
\label{ex5.5-1}
\end{align}
where $k$ is a positive parameter and $\theta = \pi \sin(x) \sin(y)$. We take $\V{b} = 0$ and $c = 0$ in (\ref{eigen-1}).

We first take $k = 10$. BAMG is able to generate meshes satisfying (\ref{thm:irreducibleM-2}) for this case.
A mesh and corresponding principal eigenfunction are shown in Fig.~\ref{ex5.5-f1}. Once again, no undershoot is
observed. The error in the computed $\lambda_1$ is shown in Fig.~\ref{ex5.5-f3}.

Next, we consider a more anisotropic case with $k = 100$. For this case, BAMG is not able to produce a mesh satisfying
the mesh condition (\ref{thm:irreducibleM-2}).  A generated mesh and corresponding principal eigenfunction are
plotted in Fig.~\ref{ex5.5-f2}. Interestingly, no undershoot is observed in this case although the stiffness matrix
is not an $M$-matrix. This indicates that the $M$-matrix requirement
(which is a sufficient requirement in Theorem~\ref{thm:fem-eigen}) can be replaced with a weaker condition.
The error in the computed $\lambda_1$ is shown in Fig.~\ref{ex5.5-f3} to have a second order convergence rate.

\end{exam}

\section{Conclusions and further comments}
\label{SEC:conclusions}

In the previous sections we have studied the P1 finite element approximation of the eigenvalue problem
of second-order elliptic differential operators subject to the Dirichlet boundary condition. The focus is on
the preservation of some basic properties of the principal eigenvalue and eigenfunctions.
It has been shown in Theorem~\ref{thm:fem-eigen} that if the stiffness matrix is an irreducible $M$-matrix,
the algebraic eigenvalue problem resulting from the P1 finite element discretization preserves most basic
properties of the principal eigenvalue and eigenfunctions of the continuous problem. These properties include 
the principal eigenvalue being real and simple and the corresponding eigenfunctions being
either positive or negative inside the physical domain. The mesh conditions leading to
such a stiffness matrix have been investigated and the main result is stated in Theorem~\ref{thm:irreducibleM}.
Roughly speaking, the theorem states that if the mesh is simplicial, acute (in 2D this condition
can be replaced by the Delaunay condition) when measured in the metric specified by the inverse
of the diffusion matrix, and interiorly connected, then the stiffness matrix is an irreducible $M$-matrix.

Numerical examples have been presented to verify the theoretical findings. They also show that
when the stiffness matrix is not an $M$-matrix, there is no guarantee that the resulting algebraic
eigenvalue problem preserve the basic properties of the principal eigenvalue and eigenfunctions.
Particularly, the eigenfunctions corresponding to the smallest eigenvalue may change sign 
and even more,  the smallest eigenvalue (in modulus) may not necessarily be real for
nonsymmetric operators. Furthermore, numerical results show that those basic properties
can be preserved for some non-$M$-matrix situations. This indicates that the $M$-matrix requirement
may be weakened. A possibility is to use generalized $M$-matrices \cite{ES2008} although it is not
obvious how conditions for generalized $M$-matrices can directly result in mesh conditions
that can be used in practical computation.
Finally, Example~\ref{ex5.5} shows that it is challenging to generate
meshes satisfying the conditions in Theorem~\ref{thm:irreducibleM}
for a general diffusion matrix. How to generate such meshes deserves more investigations in the future.

\vspace{10pt}

{\bf Acknowledgment.} The work was partially supported by the NSF under grant DMS-1115118.
The author would like to thank Weishi Liu, Erik Van Vleck, and Hongguo Xu for useful discussion
during the preparation of this work.


\begin{thebibliography}{10}

\bibitem{ABHFDV02}
D.~Ait-Ali-Yahia, G.~Baruzzi, W.~G. Habashi, M.~Fortin, J.~Dompierre, and M.-G.
  Vallet.
\newblock Anisotropic mesh adaptation: towards user-independent,
  mesh-independent and solver-independent {CFD}. {P}art {II}: {S}tructured
  grids.
\newblock {\em Int. J. Numer. Meth. Fluids}, 39:657--673, 2002.

\bibitem{BO1991}
I.~Babu{\v{s}}ka and J.~Osborn.
\newblock Eigenvalue problems.
\newblock In {\em Handbook of numerical analysis, {V}ol.\ {II}}, Handb. Numer.
  Anal., II, pages 641--787. North-Holland, Amsterdam, 1991.

\bibitem{BO1989}
I.~Babu{\v{s}}ka and J.~E. Osborn.
\newblock Finite element-{G}alerkin approximation of the eigenvalues and
  eigenvectors of selfadjoint problems.
\newblock {\em Math. Comp.}, 52:275--297, 1989.

\bibitem{BP94}
A.~Berman and R.~J. Plemmons.
\newblock {\em Nonnegative Matrices in the Mathematical Sciences}.
\newblock Society for Industrial and Applied Mathematics, Philadelphia, 1994.

\bibitem{BCER95}
M.~Bern, L.~P. Chew, D.~Eppstein, and J.~Ruppert.
\newblock Dihedral bounds for mesh generation in high dimensions.
\newblock In {\em Proc. 6th ACM-SIAM Symp. Discrete Algorithms}, pages 89--196,
  1995.

\bibitem{BEG94}
M.~Bern, D.~Eppstein, and J.~Gilbert.
\newblock Provably good mesh generation.
\newblock {\em J. Comp. System Sciences}, 48:384--409, 1994.

\bibitem{BdBSW1966}
G.~Birkhoff, C.~de~Boor, B.~Swartz, and B.~Wendroff.
\newblock Rayleigh-{R}itz approximation by piecewise cubic polynomials.
\newblock {\em SIAM J. Numer. Anal.}, 3:188--203, 1966.

\bibitem{Bof2010}
D.~Boffi.
\newblock Finite element approximation of eigenvalue problems.
\newblock {\em Acta Numer.}, 19:1--120, 2010.

\bibitem{Boffi2010b}
D.~Boffi, F.~Gardini, and L.~Gastaldi.
\newblock Some remarks on eigenvalue approximation by finite elements.
\newblock In J.~Blowey and M.~Jensen, editors, {\em Frontiers in Numerical
  Analysis Ð Durham 2010}, volume~85 of {\em Lecture Notes in Computational
  Science and Engineering}, pages 1 -- 77, Berlin, Heidelberg, 2010.
  Springer-Verlag.

\bibitem{BGHLS97}
H.~Borouchaki, P.~L. George, P.~Hecht, P.~Laug, and E.~Saletl.
\newblock {D}elaunay mesh generation governed by metric specification: {P}art
  {I}. {A}lgorithms.
\newblock {\em Fin. Elem. Anal. Des.}, 25:61--83, 1997.

\bibitem{BGM97}
H.~Borouchaki, P.~L. George, and B.~Mohammadi.
\newblock {D}elaunay mesh generation governed by metric specification: {P}art
  {II}. {A}pplications.
\newblock {\em Fin. Elem. Anal. Des.}, 25:85--109, 1997.

\bibitem{BH96}
F.~J. Bossen and P.~S. Heckbert.
\newblock A pliant method for anisotropic mesh generation.
\newblock In {\em Proceedings, 5th International Meshing Roundtable}, pages
  63--74, Sandia National Laboratories, Albuquerque, NM, 1996.
\newblock Sandia Report 96-2301.

\bibitem{BKK08}
J.~Brandts, S.~Korotov, and M.~K\v{r}\'i\v{z}ek.
\newblock The discrete maximum principle for linear simplicial finite element
  approximations of a reaction-diffusion problem.
\newblock {\em Lin. Alg. Appl.}, 429:2344--2357, 2008.

\bibitem{BKKS09}
J.~Brandts, S.~Korotov, M.~K\v{r}\'i\v{z}ek, and J.~\v{S}olc.
\newblock On nonobtuse simplicial partitions.
\newblock {\em SIAM Rev.}, 51:317--335, 2009.

\bibitem{BE04}
E.~Burman and A.~Ern.
\newblock Discrete maximum principle for {G}alerkin approximations of the
  {L}aplace operator on arbitrary meshes.
\newblock {\em C. R. Acad. Sci. Paris}, Ser.I 338:641--646, 2004.

\bibitem{CHMP97}
M.~J. Castro-D\'{i}az, F.~Hecht, B.~Mohammadi, and O.~Pironneau.
\newblock Anisotropic unstructured mesh adaption for flow simulations.
\newblock {\em Int. J. Numer. Meth. Fluids}, 25:475--491, 1997.

\bibitem{Cia70}
P.~G. Ciarlet.
\newblock Discrete maximum principle for finite difference operators.
\newblock {\em Aequationes Math.}, 4:338--352, 1970.

\bibitem{Cia78}
P.~G. Ciarlet.
\newblock {\em The Finite Element Method for Elliptic Problems}.
\newblock North-Holland, Amsterdam, 1978.

\bibitem{CR73}
P.~G. Ciarlet and P.-A. Raviart.
\newblock Maximum principle and uniform convergence for the finite element
  method.
\newblock {\em Comput. Meth. Appl. Mech. Engrg.}, 2:17--31, 1973.

\bibitem{DXZ2008}
X.~Dai, J.~Xu, and A.~Zhou.
\newblock Convergence and optimal complexity of adaptive finite element
  eigenvalue computations.
\newblock {\em Numer. Math.}, 110:313--355, 2008.

\bibitem{DZ2008}
X.~Dai and A.~Zhou.
\newblock Three-scale finite element discretizations for quantum eigenvalue
  problems.
\newblock {\em SIAM J. Numer. Anal.}, 46:295--324, 2007/08.

\bibitem{BKOS00}
M.~de~Berg, M.~van Kreveld, M.~Overmars, and O.~Schwarzkopf.
\newblock {\em Computational Geometry}.
\newblock Springer, Berlin, 2000.

\bibitem{DVBFH02}
J.~Dompierre, M.-G. Vallet, Y.~Bourgault, M.~Fortin, and W.~G. Habashi.
\newblock Anisotropic mesh adaptation: towards user-independent,
  mesh-independent and solver-independent {CFD}. {P}art {III}: {U}nstructured
  meshes.
\newblock {\em Int. J. Numer. Meth. Fluids}, 39:675--702, 2002.

\bibitem{DDS04}
A.~Dr\v{a}g\v{a}nescu, T.~F. Dupont, and L.~R. Scott.
\newblock Failure of the discrete maximum principle for an elliptic finite
  element problem.
\newblock {\em Math. Comp.}, 74:1--23, 2004.

\bibitem{ES2008}
A.~Elhashash and D.~B. Szyld.
\newblock Generalizations of {$M$}-matrices which may not have a nonnegative
  inverses.
\newblock {\em Lin. Alg. Appl.}, 429:2435 -- 2450, 2008.

\bibitem{ESU04}
D.~Eppstein, J.~M. Sullivan, and A.~Ungor.
\newblock Tiling space and slabs with acute tetrahedra.
\newblock {\em Comput. Geom. Theory \& Appl.}, 27:237--255,
  2004.

\bibitem{Eva98}
L.~C. Evans.
\newblock {\em Partial Differential Equations}.
\newblock American Mathematical Society, Providence, Rhode Island, 1998.
\newblock Graduate Studies in Mathematics, Volume 19.

\bibitem{Fix1973}
G.~J. Fix.
\newblock Eigenvalue approximation by the finite element method.
\newblock {\em Advances in Math.}, 10:300--316, 1973.

\bibitem{GS98}
R.~V. Garimella and M.~S. Shephard.
\newblock Boundary layer meshing for viscous flows in complex domain.
\newblock In {\em Proceedings, 7th International Meshing Roundtable}, pages
  107--118, Sandia National Laboratories, Albuquerque, NM, 1998.

\bibitem{GL09}
S.~G\H{u}nter and K.~Lackner.
\newblock A mixed implicit-explicit finite difference scheme for heat transport
  in magnetised plasmas.
\newblock {\em J. Comput. Phys.}, 228:282--293, 2009.

\bibitem{GYK05}
S.~G\H{u}nter, Q.~Yu, J.~Kruger, and K.~Lackner.
\newblock Modelling of heat transport in magnetised plasmas using non-aligned
  coordinates.
\newblock {\em J. Comput. Phys.}, 209:354--370, 2005.

\bibitem{HDBAFV00}
W.~G. Habashi, J.~Dompierre, Y.~Bourgault, D.~Ait-Ali-Yahia, M.~Fortin, and
  M.-G. Vallet.
\newblock Anisotropic mesh adaptation: towards user-independent,
  mesh-independent and solver-independent {CFD}. {P}art {I}: {G}eneral
  principles.
\newblock {\em Int. J. Numer. Meth. Fluids}, 32:725--744, 2000.

\bibitem{Hec97}
F.~Hecht.
\newblock {BAMG} -- {B}idimensional {A}nisotropic {M}esh {G}enerator homepage.
\newblock {http://www.ann.jussieu.fr/$\sim$hecht/ftp/bamg/}, 1997.

\bibitem{HJ1985}
R.~A. Horn and C.~A. Johnson.
\newblock {\em Matrix Analysis}.
\newblock Cambridge University Press, Cambridge, London, 1985.

\bibitem{HuJun2012}
J.~Hu, Y.~Huang, and Q.~Lin.
\newblock Lower bounds for eigenvalues of elliptic operators -- by
  nonconforming finite element methods.
\newblock 2012.
\newblock (arXiv: 1112.1145).

\bibitem{Hua01b}
W.~Huang.
\newblock Variational mesh adaptation: isotropy and equidistribution.
\newblock {\em J. Comput. Phys.}, 174:903--924, 2001.

\bibitem{Hua06}
W.~Huang.
\newblock Mathematical principles of anisotropic mesh adaptation.
\newblock {\em Comm. Comput. Phys.}, 1:276--310, 2006.

\bibitem{Hua10}
W.~Huang.
\newblock Discrete maximum principle and a delaunay-type mesh condition for
  linear finite element approximations of two-dimensional anisotropic diffusion
  problems.
\newblock {\em Numer. Math. Theory Meth. Appl.}, 4:319--334, 2011.
\newblock (arXiv:1008.0562).

\bibitem{HR11}
W.~Huang and R.~D. Russell.
\newblock {\em Adaptive Moving Mesh Methods}.
\newblock Springer, New York, 2011.
\newblock Applied Mathematical Sciences Series, Vol. 174.

\bibitem{KK09}
J.~Kar{\'a}tson and S.~Korotov.
\newblock An algebraic discrete maximum principle in {H}ilbert space with
  applications to nonlinear cooperative elliptic systems.
\newblock {\em SIAM J. Numer. Anal.}, 47:2518--2549, 2009.

\bibitem{KKK07}
J.~Kar\'atson, S.~Korotov, and M.~K\v{r}\'i\v{z}ek.
\newblock On discrete maximum principles for nonlinear elliptic problems.
\newblock {\em Math. Comput. Sim.}, 76:99--108, 2007.

\bibitem{KSS09}
D.~Kuzmin, M.~J. Shashkov, and D.~Svyatskiy.
\newblock A constrained finite element method satisfying the discrete maximum
  principle for anisotropic diffusion problems.
\newblock {\em J. Comput. Phys.}, 228:3448--3463, 2009.

\bibitem{KL95}
M.~K\v{r}\'i\v{z}ek and Q.~Lin.
\newblock On diagonal dominance of stiffness matrices in {3D}.
\newblock {\em East-West J. Numer. Math.}, 3:59--69, 1995.

\bibitem{LePot09}
C.~Le~Potier.
\newblock A nonlinear finite volume scheme satisfying maximum and minimum
  principles for diffusion operators.
\newblock {\em Int. J. Finite Vol.}, 6:20 pp., 2009.

\bibitem{Let92}
F.~W. Letniowski.
\newblock Three-dimensional {D}elaunay triangulations for finite element
  approximations to a second-order diffusion operator.
\newblock {\em SIAM J. Sci. Stat. Comput.}, 13:765--770, 1992.

\bibitem{LH10}
X.~P. Li and W.~Huang.
\newblock An anisotropic mesh adaptation method for the finite element solution
  of heterogeneous anisotropic diffusion problems.
\newblock {\em J. Comput. Phys.}, 229:8072--8094, 2010 (arXiv:1003.4530v2).

\bibitem{LSS07}
X.~P. Li, D.~Svyatskiy, and M.~Shashkov.
\newblock Mesh adaptation and discrete maximum principle for {2D} anisotropic
  diffusion problems.
\newblock Technical Report LA-UR 10-01227, Los Alamos National Laboratory, Los
  Alamos, NM, 2007.

\bibitem{LSSV07}
K.~Lipnikov, M.~Shashkov, D.~Svyatskiy, and Y.~Vassilevski.
\newblock Monotone finite volume schemes for diffusion equations on
  unstructured triangular and shape-regular polygonal meshes.
\newblock {\em J. Comput. Phys.}, 227:492--512, 2007.

\bibitem{LS08}
R.~Liska and M.~Shashkov.
\newblock Enforcing the discrete maximum principle for linear finite element
  solutions of second-order elliptic problems.
\newblock {\em Comm. Comput. Phys.}, 3:852--877, 2008.

\bibitem{LHQ2012}
C.~Lu, W.~Huang, and J.~Qiu.
\newblock Maximum principle in linear finite element approximations of
  anisotropic diffusion-convection-reaction problems.
\newblock 2012.
\newblock (arXiv:1201.3651).

\bibitem{LLX2012}
F.~Luo, Q.~Lin, and H.~Xie.
\newblock Computing the lower and upper bounds of {L}aplace eigenvalue problem:
  by combining conforming and nonconforming finite element methods.
\newblock {\em Sci. China Math.}, 55:1069--1082, 2012.

\bibitem{MM2011}
V.~Mehrmann and A.~Miedlar.
\newblock Adaptive computation of smallest eigenvalues of self-adjoint elliptic
  partial differential equations.
\newblock {\em Numer. Lin. Alg. Appl.}, 18:387--409, 2011.

\bibitem{MD06}
M.~J. Mlacnik and L.~J. Durlofsky.
\newblock Unstructured grid optimization for improved monotonicity of discrete
  solutions of elliptic equations with highly anisotropic coefficients.
\newblock {\em J. Comput. Phys.}, 216:337--361, 2006.

\bibitem{NZ2012}
A.~Naga and Z.~Zhang.
\newblock Function value recovery and its application in eigenvalue problems.
\newblock {\em SIAM J. Numer. Anal.}, 50:272--286, 2012.

\bibitem{PVMZ97}
J.~Peraire, M.~Vahdati, K.~Morgan, and O.~C. Zienkiewicz.
\newblock Adaptive remeshing for compressible flow computations.
\newblock {\em J. Comput. Phy.}, 72:449 -- 466, 1997.

\bibitem{RLSF04}
J.-F. Remacle, X.~Li, M.~S. Shephard, and J.~E. Flaherty.
\newblock Anisotropic adaptive simulation of transient flows using
  discontinuous {G}alerkin methods.
\newblock {\em Int. J. Numer. Methods Engrg.}, 62:899--923, 2005.

\bibitem{SH07}
P.~Sharma and G.~W. Hammett.
\newblock Preserving monotonicity in anisotropic diffusion.
\newblock {\em J. Comput. Phys.}, 227:123--142, 2007.

\bibitem{ShYu2011}
Z.~Sheng and G.~Yuan.
\newblock The finite volume scheme preserving extremum principle for diffusion
  equations on polygonal meshes.
\newblock {\em J. Comput. Phys.}, 230:2588--2604, 2011.

\bibitem{Sto86}
G.~Stoyan.
\newblock On maximum principles for monotone matrices.
\newblock {\em Lin. Alg. Appl.}, 78:147--161, 1986.

\bibitem{SF73}
G.~Strang and G.~J. Fix.
\newblock {\em An Analysis of the Finite Element Method}.
\newblock Prentice Hall, Englewood Cliffs, NJ, 1973.

\bibitem{WaZh11}
J.~Wang and R.~Zhang.
\newblock Maximum principle for {P1}-conforming finite element approximations
  of quasi-linear second order elliptic equations.
\newblock {\em SIAM J. Numer. Anal.}, 50:626--642, 2012.
\newblock (arXiv:1105.1466).

\bibitem{XZ2001}
J.~Xu and A.~Zhou.
\newblock A two-grid discretization scheme for eigenvalue problems.
\newblock {\em Math. Comp.}, 70:17--25, 2001.

\bibitem{XZ99}
J.~Xu and L.~Zikatanov.
\newblock A monotone finite element scheme for convection-diffusion equations.
\newblock {\em Math. Comput.}, 69:1429--1446, 1999.

\bibitem{YS00}
S.~Yamakawa and K.~Shimada.
\newblock High quality anisotropic tetrahedral mesh generation via ellipsoidal
  bubble packing.
\newblock In {\em Proceedings, 9th International Meshing Roundtable}, pages
  263--273, Sandia National Laboratories, Albuquerque, NM, 2000.
\newblock Sandia Report 2000-2207.

\bibitem{YB2011}
Y.~Yang and H.~Bi.
\newblock Two-grid finite element discretization schemes based on
  shifted-inverse power method for elliptic eigenvalue problems.
\newblock {\em SIAM J. Numer. Anal.}, 49:1602--1624, 2011.

\bibitem{YZL2010}
Y.~Yang, Z.~Zhang, and F.~Lin.
\newblock Eigenvalue approximation from below using non-conforming finite
  elements.
\newblock {\em Sci. China Math.}, 53:137--150, 2010.

\bibitem{YuSh2008}
G.~Yuan and Z.~Sheng.
\newblock Monotone finite volume schemes for diffusion equations on polygonal
  meshes.
\newblock {\em J. Comput. Phys.}, 227:6288--6312, 2008.

\bibitem{ZZS2013}
Y.~Zhang, X.~Zhang, and C.-W. Shu.
\newblock Maximum-principle-satisfying second order discontinuous {G}alerkin
  schemes for convection-diffusion equations on triangular meshes.
\newblock {\em J. Comput. Phys.}, 234:295 -- 316, 2013.

\end{thebibliography}

\end{document}